\documentclass[11pt]{article}
\usepackage{graphicx}
\usepackage{amsmath}
\usepackage{amsfonts}
\usepackage{epsfig}

\newcommand{\be}{\begin{equation}}
\newcommand{\ee}{\end{equation}}
\newcommand{\beqn}{\begin{eqnarray}}
\newcommand{\eeqn}{\end{eqnarray}}
\newcommand{\beqns}{\begin{eqnarray*}}
\newcommand{\eeqns}{\end{eqnarray*}}
\newcommand{\h}{\widehat}

\newcommand{\fr}[1]{(\ref{#1})}
\newcommand{\supp}{\mbox{supp}\ }
\newcommand{\lkr}{\left(}
\newcommand{\rkr}{\right)}
\newcommand{\lkv}{\left[}
\newcommand{\rkv}{\right]}
\newcommand{\lfi}{\left\{}
\newcommand{\rfi}{\right\}}
\newcommand{\EB}{EB\ }

\newcommand{\EE}{\ensuremath{{\mathbb E}}}
\newcommand{\II}{\ensuremath{{\mathbb I}}}
\newcommand{\PP}{\ensuremath{{\mathbb P}}}

\newtheorem{example}{Example}
\newtheorem{lemma}{Lemma}
\newtheorem{theorem}{Theorem}

\newtheorem{remark}{Remark}

\title{\Large{\bf Adaptive  Nonparametric Empirical Bayes Estimation 
Via Wavelet Series: the  Minimax Study  }}

\author{
\large{ Rida Benhaddou and Marianna Pensky}  
  \\ \\
Department of Mathematics, University of Central Florida   } 

\date{}

\begin{document}

\maketitle

\begin{abstract}
In the present paper, we derive  lower bounds for the risk of the nonparametric 
empirical Bayes estimators.   In order to attain the optimal convergence rate, we
propose  generalization of the linear  empirical Bayes  estimation method
which takes advantage of  the flexibility of the wavelet techniques.
We present an empirical Bayes   estimator as a wavelet series  expansion  
and estimate coefficients by minimizing the prior risk of the estimator.
As a result,   estimation of wavelet coefficients requires solution 
of a well-posed low-dimensional sparse system of linear equations. 
The dimension of the system depends on the size of wavelet support  
and smoothness of the Bayes estimator.  An adaptive choice of the resolution level 
is carried out using Lepski~(1997) method. The method is computationally efficient and provides 
asymptotically optimal adaptive EB estimators. 
The theory is supplemented by numerous examples.\\

\end{abstract}

\section{Introduction}
\label{sec:introduction}
\setcounter{equation}{0}

    Empirical Bayes (EB)  methods are   estimation techniques in which the prior distribution, 
in the standard Bayesian sense, is estimated from the data. They are powerful tools, in particular, 
when data are generated by repeated execution of the same type of experiment. The EB methods are 
directly related to the standard Bayes models but there is difference in perspective between the 
two in the sense that, in the standard Bayesian approach, the prior distribution 
is assumed to be fixed before any data are observed, whereas, in the EB setting, the prior 
distribution, in some way or another,  is estimated from the observed data.

Consider the following setting. One observes independent two-dimensional  random  vectors
$(X_1,  \theta_1)$  $,..., (X_n, \theta_n),$
where  each $\theta_i$ is distributed
according to some unknown prior pdf  $g$ and, given $\theta_i =\theta$,   observation 
$X_i$ has the known   conditional density function
$q(x \mid\theta)$, so that, each pair $(X_i, \theta_i)$ has an absolutely continuous 
distribution with the density function $q(x|\theta) g(\theta)$.  In each pair, the first component
is observable, but the second is not.
After  the $(n+1)$-th observation $y \equiv X_{n+1}$ is taken,
the goal is to estimate  $t \equiv \theta_{n+1}$.

If the prior density $g(\theta)$ were known, then the   Bayes estimator of $\theta_{n+1}$ which 
delivers the minimal mean squared risk would be given by the following equation
 \beqn  \label{pm}
  t(y)=\frac{\int^{\infty}_{-\infty} \theta q(y\mid\theta){g}(\theta)d\theta}
{\int ^{\infty}_{-\infty} q(y\mid\theta){g}(\theta)d\theta}.   
\eeqn
(see e.g., Carlin and Louis~(2001), or Maritz and Lwin~(1989)).\\

 Since the prior density $g(\theta)$ is unknown, an EB estimator $\h{t}(y; X_1, X_2, \cdots, X_n)$ 
has to be used. Using notations 
 \beqn 
p(y) & = & \int_{-\infty}^{\infty}   q(y|\theta)  g(\theta) d\theta, \label{py}\\
\Psi(y) & = & \int_{-\infty}^{\infty}  \theta q(y|\theta) g(\theta) d\theta. \label{Phiy}
\eeqn
 $t(y)$ can be rewritten as 
\be \label{eb1}
t(y)=  \Psi(y)\big/ p(y).
\ee

There is a variety of methods which allow to estimate $t(y)$ on the basis of
observations $y; X_1, \cdots, X_n$.  
After Robbins (1955, 1964) formulated   EB estimation approach,  many statisticians have been working on 
developing  EB methods. The comprehensive list of references as well as numerous examples of applications of EB techniques can 
be found in  Carlin and  Louis  (2000) or Maritz  and Lwin (1989).

The EB techniques can be divided into two groups: parametric and nonparametric. 
Parametric \EB methods require that the parametric form  of a family, to which the 
prior distribution belongs, is specified a priori. 
The past data is then used to estimate the values of the unknown parameters, usually using 
the maximum likelihood  approach (see, e.g., Louis~(1984), Morris~(1983), Casella~(1985) 
and Efron and Morris~(1977) among others).

In nonparametric \EB estimation,  prior distribution is completely unspecified. 
One of the approaches to nonparametric EB estimation is based on estimation of the numerator 
and the denominator in the ratio  in  \fr{eb1}. This approach was introduced by 
Robbins  (1955, 1964) himself   and later developed by a number of authors 
(see, e.g.,    Brown and Greenshtein (2009), Datta  (1991, 2000), 
Ma and Balakrishnan (2000), Nogami  (1988),  Pensky  (1997a,b),  
Raykar and Zhao (2011),  Singh  (1976, 1979) and Walter and Hamedani  (1991) 
among others). The method provides  estimators with good convergence rates, 
however, it requires  relatively tedious three-step procedure: estimation 
of the top and the bottom of the fraction and then the fraction itself.
In particular, one of the approaches is to take advantage of the fact 
that, in the case of a one-parameter exponential family, 
the numerator can be expressed as the derivative of the denominator.
Wavelets provide  an opportunity to construct adaptive wavelet--based EB estimators   
with better computational properties in this framework
(see, e.g., Huang  (1997)  and Pensky (1998, 2000, 2002))   but 
the necessity of estimation of the ratio in \fr{eb1} remains. 
Another nonparametric approach developed in  Jiang and Zhang (2009),
is based on application of nonparametric MLE technique which is computationally  extremely
demanding.

In 1983, Robbins introduced a much more simple, local nonparametric EB method, 
 linear EB estimation. Robbins (1983) suggested to approximate  Bayes estimator $t(y)$ 
locally by  a linear function of $y$ and to determine the coefficients of $t(y)$
by  minimizing the expected squared difference between $t(y)$ and $\theta$, 
with subsequent estimation of the coefficients on the basis of observations $X_1,...,X_n$.
The technique is extremely efficient computationally and  was immediately put to practical use, 
for instance, for prediction of the finite population mean (see, e.g., Ghosh and Meeden (1986),
Ghosh and Lahiri  (1987) and  Karunamuni and Zhang (2003)).

However, a linear EB  estimator   has a   large bias since, due to its very simple form, 
it has a limited ability to approximate the Bayes  estimator $t(y)$. For this reason, 
linear EB estimators are optimal only in the class of estimators {\it linear} in $y$. 
To overcome this defect, Pensky and Ni (2000)  extended approach of Robbins (1983)  
to approximation of $t(y)$ by  algebraic polynomials.  However, although  the 
polynomial-based EB estimation provides significant improvement in the convergence rates in comparison with 
the linear EB estimator, the system of linear equations resulting from the method is badly 
conditioned which leads to computational difficulties and loss of precision.

To overcome those difficulties, Pensky and Alotaibi (2005) proposed  
to replace polynomial approximation of the Bayes estimator $t(y)$ by 
its approximation via wavelets, in particular, by expansion over scaling 
functions at the resolution level $m$. The method  exploits de-correlating property 
of wavelets and leads to a low-dimensional well-posed sparse system  of linear equations.
The paper also treated the issue of locally   optimal choice of resolution level as  $n \rightarrow \infty$: 
if the resolution level is selected correctly, in accordance with the smoothness of the Bayes estimator,
 then the suggested EB estimator attains the best convergence rates which can be obtained by
application of wavelet-based \EB\ estimator.

However, smoothness of the Bayes estimator  $t(y)$ is hard to assess.
For this reason, the EB estimator of Pensky and Alotaibi (2005)  is non-adaptive.
One of the possible ways of achieving adaptivity would be to replace the 
linear scaling function based approximation by a traditional wavelet expansion with 
subsequent thresholding of wavelet coefficients. The deficiency of this approach, however, is 
that it yields  the system of equations which is much less sparse and is growing in size 
with the number of observations $n$.

The present paper has two main objectives. The first one is to derive lower bounds for the posterior risk
of a nonparametric empirical Bayes estimator. In spite of a fifty years  long history of empirical Bayes estimation methods, 
a general lower  bound  for the risk of an empirical Bayes estimators has not been derived so far.
Only some particular cases of the problem were investigated.
Specifically, Penskaya (1995) obtain lower bounds for the posterior risk of nonparametric empirical Bayes estimators 
of a location parameter. Li, Gupta and  Liese (2005) obtained lower bounds for the risk of empirical Bayes estimators 
in the exponential families. However, since their  lower bound  was of the form $C n^{-1}$, practically no estimator could attain it.  
Construction of the lower bounds for the risk was attempted also in the empirical Bayes  two-action problem in the case of 
a continuous one-parameter exponential family. Karunamuni (1996) published the paper on the subject but his results were proved to be inaccurate,
at least in the case of the normal distirbution, when Liang (2000) constructed an estimator with the convergence rates below the lower bound
for the risk. Pensky  (2003) derived lower bounds for the loss in the  empirical Bayes two-action problem involving normal means.
Nevertheless, no general theory for construction of the lower bounds for the prior or posterior  risk of an empirical Bayes estimator
  has ever been developed so far. In what follows, we construct lower bounds for the posterior risk of an empirical Bayes estimator 
under a general assumption that the marginal density $p(x)$ given by formula \fr{py} is continuously differentiable in the neighborhood of $y$.

The second purpose of this paper is to provide an adaptive version of the wavelet EB estimator 
developed in Pensky and Alotaibi (2005). In particular, we preserve the linear structure of the estimator. 
However, since expansion over scaling functions at the resolution level $m$  leads to excessive variance when resolution level $m$ is too high 
and disproportionately large bias when $m$ is too  small, we choose the resolution level using   
methodology introduced by  Lepski  (1991) and further developed by  Lepski, Mammen and Spokony (1997). 
The resulting estimator is adaptive and attains optimal convergence rates (within   a logarithmic 
factor of $n$). In addition, it has an advantage of computational efficiency since it is based 
on the solution of low-dimensional sparse system of  linear equations the matrix of which 
tends to a scalar multiple of an identity matrix as the scale $m$ grows. 
The theory is supplemented by numerous examples that demonstrate how the estimator can be 
implemented for various types of distribution families.

The rest of the paper is organized as follows.  
Section~\ref{sec:estimation} introduces EB estimation algorithm. 
Section~\ref{sec:conv_rates} assesses  estimation error and describes the choice 
of the resolution level which delivers the best possible  convergence rates when the degree of smoothness 
of Bayes estimator $t(y)$ is known. Section~\ref{sec:lower_bounds} derives minimax lower bounds 
for the posterior risk, so that we can verify that the EB estimators constructed in the paper
are indeed asymptotically optimal as $n \to \infty$. Section~\ref{sec:Lepskii_method} 
discusses adaptive choice of the resolution level using Lepski method and proves asymptotic 
optimality of the resulting EB estimator which is based on this choice. Section~\ref{sec:examples} provides examples 
of construction of EB estimators for a variety of distribution families. Section~\ref{sec:discussion}
concludes the paper with discussion. Finally, Section~\ref{sec:proofs} contains the proofs 
of the statements in the paper.

\section {EB estimation algorithm }
\label{sec:estimation}
\setcounter{equation}{0}

In order to construct an estimator of $t(y)$ defined in \fr{eb1},
choose  a twice continuously differentiable scaling function $ \varphi$  
with bounded support and $s$ vanishing moments, so that
\beqn
 \supp \varphi& \in & [M_1, M_2],
\label{wavco2} \\
\int_{-\infty}^{\infty} x^{\aleph}  \sum_{k \in Z}  \varphi(x-k) \varphi(z-k)   dx & = & z^{\aleph}, \ \ 0   \leq   \aleph \leq s-1. 
\label{wavco1}  
\eeqn 
(see e.g., Walter and Shen~(2001)).\\

Approximate $t(y)$ by a wavelet series, for some fixed $m \geq 0$, 
\be
t_m(y) = \sum_{k \in Z}  a_{m,k}\, \varphi_{m,k} (y) 
\label{tmy}
\ee
where $\varphi_{m,k} (y) = 2^{m/2} \varphi(2^m y -k)$, and  estimate coefficients of $t_m(y) $ by 
minimizing the integrated mean squared error 
\be
\min_{a_{m,k}} \left \{ \int_{-\infty}^{\infty} \int_{-\infty}^{\infty} \lkv \sum_{k \in Z}  a_{m,k} \varphi_{m,k} (y) - z \rkv^2 q(y|z) g(z) dz\, dy \right\}.
\label{mindi}
\ee  
Taking derivatives of the last expression with respect to $a_{m,j}$ 
and equating them to zero, we obtain the system of linear equations 
\be \label{linsys}
B_m a_m = c_m 
\ee 
with 
\beqn
(B_m)_{j,k} = B_{j,k} & = & \int_{-\infty}^{\infty} \varphi_{m,k} (x) \varphi_{m,j} (x) p(x) dx = 
\EE\lkv \varphi_{m,k} (X) \varphi_{m,j} (X) \rkv,
\label{bjk}\\
(c_m)_j = c_j  & = & \int_{-\infty}^{\infty} \varphi_{m,j} (x) \Psi (x) dx. 
\label{cj}
\eeqn
where   we use the symbol $\EE$ for expectation over 
the joint distribution of $X_1, X_2, \cdots, X_n$.  The expectations over 
all other distributions are represented in the  integral forms.
Also, in what follows, we suppress index $m$ in notations of matrix $B_m =B$ and vector $c_m =c$
unless this leads to a confusion.

System \fr{linsys} is an infinite system of equations. However, 
since we are interested in estimating $t(x)$   locally at $x=y$, 
we shall keep  only indices $k,j \in K_{m,y}$ such that 
\be
K_{m,y} = \lfi k \in Z:\ 2^m y - M_2 - s(M_2 - M_1) \leq k \leq 2^m y - M_1 + s(M_2 - M_1) \rfi,
\label{indset}
\ee
where $s$ is the number of vanishing moments of   scaling function $\varphi$ (see \fr{wavco1}).
Observe that expansion \fr{tmy} actually contains only coefficients 
$a_{m,k}$ with $2^m y - M_2 \leq k \leq 2^m y - M_1$,
however, in order to ensure fast convergence of the bias of the estimator to zero,
 we need to keep more terms in the system of equations \fr{linsys}  (see Lemma A.3 in Pensky and Alotaibi~(2005) for more details).

The entries \fr{bjk} of the matrix $B$ are unknown and can be estimated by 
sample means 
\be
\h{B}_{j,k} = n^{-1}\displaystyle\sum\limits^n_{l=1}\lkv \varphi_{m,k}(X_l) \varphi_{m,j} (X_l) \rkv.
\label{hbjk}
\ee
In order to estimate $c_j$, one needs to find functions $u_{m,j} (x)$ such that for any $\theta$
\be
\int_{-\infty}^{\infty} q(x|\theta) u_{m,j} (x) dx = \int_{-\infty}^{\infty} \theta q(x|\theta) \varphi_{m,j} (x) dx.
\label{eq1}
\ee
Then, multiplying both sides of \fr{eq1} by $g(\theta)$ and 
integrating over $\theta$, we obtain
$$
\EE  u_{mj} (X) =  \int_{-\infty}^{\infty} u_{mj} (x) p(x) dx = \int_{-\infty}^{\infty} \varphi_{m,j} (x) \Psi (x) dx = c_j.
$$
Note that functions $u_{m,j} (x)$ are the same functions that appeared 
in the wavelet estimator of the numerator $\Psi (y)$ of the EB estimator \fr{eb1}
in  Pensky (1997a, 1998) where it was demonstrated that construction of $u_{m,j} (x)$ is
possible in many particular cases. In Section \ref{sec:examples} we show that 
solutions of equation \fr{eq1} can be easily obtained for a variety of distribution families.

Once functions $u_{m,j} (x)$ are derived, coefficients $c_j$ can be estimated by 
\be
\h{c_j} = n^{-1} \displaystyle\sum\limits^n_{l=1} u_{m,j} (X_l) 
\label{hcj}
\ee
and system \fr{linsys} is replaced by $\h{B} \h{a} = \h{c}$.
However, though estimators $\hat{B}$ and $\hat{c}$ converge in mean squared sense  
to $B$ and $c$, respectively, the estimator $\h{a} = \h{B}^{-1}\, \h{c}$ may not even have 
finite expectation. To understand this fact, note that both $\h{B}$ and $\h{c}$ are 
asymptotically normal. In one dimensional case, the ratio of two normal random variables 
has Cauchy distribution and, hence, does not have finite mean. In multivariate case the difficulty remains. 
To ensure that the estimator of $a$ has finite expectation, we choose $\delta= \delta_n >0$ and 
construct an estimator of $a$ of the form 
\be  \label{solu}
\h{a}_\delta = (\h{B}  + \delta I)^{-1}\, \h{c} 
\ee
where $I$ is the identity matrix.
Observe that matrix $\h{B}$ is nonnegative definite, so that $(\h{B}  + \delta I)$ is 
a positive definite matrix and, hence, is nonsingular. 
Solution $\h{a}_\delta$ is used for construction of the EB estimator
\be  \label{ebes}
\h{t}_m (y) = \displaystyle\sum\limits_{k \in K_{m,y}}( \h{a}_\delta)_{m, k} \, \varphi_{m,k} (y).
\ee

\section{Estimation error and convergence rates}  
\label{sec:conv_rates}
\setcounter{equation}{0}

\subsection{The posterior and the prior risks}

An EB\  estimator $\h{t}(y)$ can be characterized by its posterior risk
$$ 
R(y;  \hat{t} ) = (p(y))^{-1} \EE\int_{-\infty}^{\infty} (\hat{t}  (y) -  \theta)^2 q(y|\theta) g(\theta) d\theta 
$$  
which can be partitioned  into two components.
The first component of this sum is
$$ R(y; t(y))= \inf_f R(y; f(y))=
 (p(y))^{-1} \int_{-\infty}^{\infty} (t(y) - \theta)^2 q(y|\theta) g(\theta) d\theta, $$ 
which is independent of  $\hat{t} (y)$ and represents 
the posterior risk of the Bayes estimator (\ref{pm}).
 Thus,  we shall measure precision of an EB\  estimator \fr{ebes} by the second component
\be
R_n(y) = \EE (\hat{t}_m (y) - t(y))^2   
\label{rfun}
\ee
which represents the local error of the EB estimator at the point of observation $y$.

It must be noted that often the precision of an EB\  estimator is
described by 
$$
\EE R _n (y) = \int_{-\infty}^{\infty} \hat{R}_n (y) p(y) dy,
$$
which is the difference between the
prior risk 
$$
\EE\int_{-\infty}^{\infty} R(y; \hat{t} (y)) p(y) dy
$$ 
of the EB\  estimator $ \hat{t} (y)$ and the prior risk 
$$
\int_{-\infty}^{\infty} R(y; t(y)) p(y) dy = \inf_f \int_{-\infty}^{\infty} R(y;f(y)) p(y) dy
$$ 
of the  corresponding Bayes estimator $t(y)$.  However, the risk function \fr{rfun} 
has several advantages compared with  $\EE\hat{R}_n (y)$.

First, $R_n (y)$ enables one  to calculate the  mean squared error
for the  given  observation $y$ which is the  quantity of interest.
Note that the wavelet series  \fr{ebes} is local in a sense 
that coefficients $(\h{a}_\delta)_{m, k}$ change whenever $y$ changes,
hence, working with a local measure of the risk makes much more sense. 
Using the prior risk for the estimator which is local in nature 
prevents one from seeing advantages of this estimator. 
Second, by using the risk function \fr{rfun} we eliminate  the   influence
on the risk function of the observations having very  low   probabilities.
So, the use of   $R_n (y)$ provides a way of  getting EB\ 
estimators  with better convergence rates.
Third, posterior risk allows one to assess optimality of \EB estimators for majority of 
familiar distribution families via comparison of the convergence rate of the estimator 
with  the lower bounds for the risk derived in Pensky(1997).
 Finally, one can pursue evaluation  of the prior risk for the estimator \fr{ebes}. 
 The derivation will require assumptions similar to the ones in  Pensky  (1998)  
 and can be accomplished by standard methods.

The error \fr{rfun} is dominated by the sum of two components  
\be
R_n(y) \leq 2(R_{1n}(y)  + R_{2n}(y))
\label{part}
\ee
where the first component $R_{1n}=R_{1n}(y)$ is due to replacement of the 
Bayes estimator $t(y)$ by its wavelet representation \fr{tmy}, 
while $R_{2n}=R_{2n}(y)$ is due to replacement of vector $a= B^{-1} c$ by $\h{a}_\delta$ given by \fr{solu} 
\beqn
R_{1n} (y) & = & (t_m (y) - t(y))^2, \label{bias}\\
R_{2n} (y) & = & \EE\lkv  \displaystyle\sum\limits_{k \in K_{m,y}} [(\h{a}_{\delta})_{m, k} - a_{m,k}] \varphi_{m,k} (y) \rkv^2. 
\label{variance}
\eeqn
We shall refer to $R_{1n} (y)$ and $R_{2n} (y)$ as   the systematic  
and the random error components, respectively. Since in this paper we are using 
the posterior risk as a measure of precision of an \EB\ estimator, from now on we treat $y$ as a fixed quantity.

\subsection{The systematic error component  }

For evaluation of the  systematic error component  $R_{1n}$, let us introduce 
matrix $U_h$ and   vector  $D_h$  with components
\beqn
(U_h)_{k,l} & = & \int_{-\infty}^{\infty} z^h \varphi(z+ 2^m y-k) \varphi(z+ 2^m y-l) dz, \label{uhkl} \\
(D_h)_k & = &   \int_{-\infty}^{\infty} z^h \varphi(z+ 2^m y-k)  dz. \label{dhk} 
\eeqn
Observe that $U_h$ and $D_h$ are independent of unknown functions $p(x)$ and $\Psi(x)$, 
and that $U_0 =I$ where $I$ is the identity matrix. Denote 
\be
\Omega_{m,y} = 
\lfi x:\ |x-y| \leq 2^{-m} s (M_2 - M_1)  \rfi.
\label{ommy}
\ee
Then, the following statement is valid.

\begin{lemma}  \label{lem:syserror}
Let functions $p(x)$ and $\Psi(x)$ be $r \leq s-1$ times continuously differentiable in the neighborhood 
$\Omega_{y}$ of $y$ and let $\Omega_{m,y} \subseteq \Omega_{y}$, with $\Omega_{m,y}$ defined in \fr{ommy}. 
Then, for $R_{1n}$ defined in \fr{bias},   as $m \rightarrow \infty$,   
\be \label{R1_bound}
R_{1n} = (t_m (y) - t(y))^2 = o(2^{-2mr}).
\ee 
\end{lemma}

Let us now give some insight into the proof of the lemma. Note that,
\beqns
B_{j,k} &=& 2^m \int_{-\infty}^{\infty} \varphi(2^m x-k) \varphi(2^m x-j) p(x) dx  \\
     &=& \int_{-\infty}^{\infty} \varphi(z + 2^m y -k) \varphi(z + 2^m y -j) p(y + 2^{-m} z) dz   \\
     &=& {\sum_{h=0}^{r}} 2^{-mh} (h!)^{-1} p^{(h)} (y) U_h +  o(2^{-mr}), \nonumber
\eeqns
where $U_0 = I$ is the identity matrix. Deriving a similar representation for $c_k$, we obtain asymptotic expansions 
of matrix $B$ and vector $c$ via matrices $U_h$ and vectors $D_h$, respectively, as $m \rightarrow \infty$   
\beqn
B &=& p(y) I +  {\sum_{h=1}^{r}} 2^{-mh} (h!)^{-1} p^{(h)} (y) U_h +  o(2^{-mr}),
\label{Bexpan}\\
c & = &  2^{-m/2}\  \sum_{l=0}^{r}  2^{-ml}(l!)^{-1} \Psi^{(l)} (y) D_l +  o(2^{-mr}).
\label{cexpan}
\eeqn
Formula \fr{Bexpan} establishes that, for large values of $m$, matrix $B$ is close to $p(y) I$, so the system 
of equations \fr{linsys} is well-conditioned. Therefore, the inverse matrix is close to $p^{-1} (y) I$,
or, more precisely (see Pensky and Alotaibi (2005)), 
\be \label{B_inverse}
B^{-1} = p^{-1}(y) I - 2^{-m} p'(y)  p^{-2}(y) U_1 + o(2^{-m}).
\ee
Furthermore, if $m \rightarrow \infty$,  vector $a$ in 
\fr{linsys} tends to $2^{-m/2} [\Psi(y)/p(y)] D_0$  where  
$$
2^{-m/2} \sum_k (D_0)_k \varphi_{m,k} (y) =1
$$
for any $y$. The latter implies that the systematic error $R_{1n}$  tends to zero, as $m \rightarrow \infty$, 
at a rate $o \lkr 2^{-2mr} \rkr$.

The exact proof of Lemma \ref{lem:syserror}  can be obtained by a slight modification of the proof 
of Lemma~3.1 in Pensky and   Alotaibi (2005).
\\

\subsection{The random error component}

In order to calculate the random error component $R_{2n} (y)$ given by \fr{variance},  
introduce vectors   $\gamma^{(\varrho)} (m),\  \varrho=1,2, 3, 4.$, with components 
\be \label{gavec}
\gamma^{(\varrho)}_k (m) =  \lkv \int_{-\infty}^{\infty} u_{m,k}^{2\varrho} (x)  dx \rkv^{1/2} ,\ \ k \in K_{m,y}, \  \varrho=1,2, 3, 4.
\ee 
where $u_{m,k} (x)$ are defined in \fr{eq1}, and denote
\be \label{notation}
\gamma_m  = \|\gamma^{(1)} (m) \| 
\ee
where $\|z \|$ is the Euclidean norm of the vector $z$.
The following lemma provides an asymptotic expression for the random error component as $m,n \to \infty$.

\begin{lemma} \label{lem:randerror}
Let $\delta^2 \sim n^{-1} 2^m$. Then, under the assumptions of   Lemma  1, as $m,n \rightarrow \infty$,
the random error component $R_{2n}$ defined in  \fr{variance} is such that
\be \label{delta2}
R_{2n} = O \lkr 2^m n^{-1} (1 + \gamma_m^2) \rkr, \ \ m,n \to \infty,
\ee
provided   $2^m n^{-1} \rightarrow 0$ and $\|\gamma^{(2)} (m) \|^2\, 2^{2m} = o(n^3)$
as $n \rightarrow \infty$.   
\end{lemma}

Proofs of this and other statements in the paper are given in Section \ref{sec:proofs}
Proofs.\\

Observe that the values of $\gamma_k^{(\varrho)} (m)$ are independent of the unknown density $g(\theta)$ and can be calculated 
explicitly. In Section \ref{sec:examples}, we  bring examples of construction of functions $u_{m,k} (x)$ as well as 
the asymptotic expressions for $\gamma_k^{(\varrho)} (m)$, $\varrho=1,2$, for some common special cases
(location parameter  family, scale parameter family, one-parameter exponential family). 
In vast majority of situations,  including the ones studied in Section \ref{sec:examples}, 
$\gamma_ m^2$ is  bounded above  by the following   expression 
\be \label{gamma_param_form}
\gamma_m^2 \leq C_{\gamma} 2^{\alpha m}  \ \  C_{\gamma} >0, \ \alpha \in \mathbb R,
\ee
where $\alpha$ and $C_{\gamma}$ are the absolute constants independent of $m$.

It follows from Lemmas \ref{lem:syserror} and \ref{lem:randerror} that $R_n$ given in \fr{part}
is minimized at the optimal resolution level   $m=m_0$ where both errors are balanced
\be \label{mnopt}
m_0 = \arg \min \lkr  n^{-1} 2^m (\gamma_m^2 +1)  + 2^{-2mr} \rkr.
\ee
In particular, under assumption \fr{gamma_param_form}, as $n \rightarrow \infty$, 
$m_0$ is such that
\be \label{m_opt}
2^{m_0} \asymp n^{\frac{1}{2r + 1+ \max(\alpha,0) }}.
\ee
Here, we denote  $a_n \asymp b_n$ for two sequences,  $\{ a_n \}$ and  $\{ b_n \}$, $n = 1,2, \cdots,$ of positive real numbers 
if there exist     $C_1$ and $C_2$  independent of $n$ 
such that $0 < C_1 < C_2 < \infty$ and $C_1 \leq \alpha_n/\beta_n  \leq C_2$.

Then the following statement is true.

\begin{theorem} \label{th:nonadaptive}
Let twice continuously differentiable scaling function $\varphi$ satify  \fr{wavco2} and \fr{wavco1}.
Let functions $p(x)$ and $\Psi(x)$ be $r$   
times continuously differentiable in the neighborhood $\Omega_{y}$ of $y$ such that
$\Omega_{m,y} \subseteq \Omega_{y}$ where $\Omega_{m,y}$ is defined in \fr{ommy} and   $1/2 \leq r \leq (s-1)$.
Choose $m_0$ according to \fr{mnopt} and let $\delta$ in \fr{solu} be such that $\delta^2 \sim n^{-1} 2^{m_0}$.
Then, for any $y$ such that $p(y)>0$,  as $n \rightarrow \infty$, 
$R_n(y)$ defined in \fr{rfun} satisfies the following  asymptotic relation
\be \label{finerror}
R_n(y) = \EE (\hat{t}_m (y) - t(y))^2  =  O \lkr 2^{-2 m_0 r} \rkr,\ \ m \to \infty, 
\ee
provided $2^m n^{-1} \to 0$ and $\|\gamma^{(2)} (m_0) \|^2 2^{2m_0} = o(n^3)$ as $n \rightarrow \infty$.
In particular, if assumption \fr{gamma_param_form} holds, then, as $n \rightarrow \infty$, one has
\be \label{asymp_error}
R_n(y) = O \lkr n^{-\frac{2r}{2r + 1+  \max(\alpha,0)}} \rkr.
\ee
\end{theorem}

\begin{remark}  \label{remark1}
{\rm In general, it is possible that instead of inequality \fr{gamma_param_form} one has 
$$ 
\gamma_m^2 \leq C_{\gamma} 2^{\alpha m} \exp \lkr b 2^{\beta m} \rkr, \quad   b,\beta, C_\gamma \geq 0.
$$
In this case, $2^{m_0} \asymp  ((2b)^{-1} \log n)^{1/\beta}$  and 
$$
R_n(y) = O \lkr  (\log n)^{-\frac{2r}{\beta}} \rkr.
$$
Note that, unlike for the case of $b=0$,  the value of $m_0$ is independent of the unknown smoothness parameter
 $r$, so the estimator is {\it adaptive}. However, the rate of convergence is very slow (logarithmic in $n$).
Fortunately, this case is not common in practice. It does not take place for any of the examples 
considered in this paper. Moreover, we failed to find an example when this situation takes place.
}
\end{remark}

\section{Minimax lower bounds for the posterior risk}  
\label{sec:lower_bounds}
\setcounter{equation}{0}

Let $y$ be a fixed point. Consider an $r$-times continuously differentiable pdf $p_0(x)$, and an $r$-times continuously 
differentiable  kernel $k(\cdot)$  with $\supp k = ( -1, 1 )$ and such that $\int k(z)dz = 0$. 
Let $p_0(\cdot)$ and $k(\cdot)$ satisfy the following assumptions:
\\

\noindent
{\bf Assumption A1}\label{assu1}:\ \  There exists $g_0(\theta)$ such that, for any $x$,
$$ 
p_0(x) = \int_{-\infty}^\infty  q(x \mid \theta ) g_0(\theta)d\theta. 
$$

\noindent
{\bf Assumption A2}\label{assu2}:\ \  There exists a function $\psi_{h, y}(\theta)$ such that, for any $x$   and $h>0$,
\be \label{psi-equation} 
k\left( \frac{x-y}{h} \right) = \int_{-\infty}^\infty q( x \mid \theta ) \psi _{h, y }(\theta) d\theta. 
\ee

\noindent
{\bf Assumption A3}\label{assu3}:\ \  Density $p_0(x)$ is such that, for any $x$, such that 
$\left| x-y \right| \leq h$ and some $C_p >0$  
$$ 
p_0(x) > 2 C_p \| k \|_{\infty}. 
$$
\\

\noindent
Denote
\begin{eqnarray}
\Psi_0(x) &=& \int_{-\infty}^\infty \theta q( x \mid \theta ) g_0(\theta)d\theta, \quad
w_{h, y}(x)   =  \int_{-\infty}^\infty \theta q(x\mid \theta) \psi_{h, y}(\theta)d\theta, 
\label{w}\\
\rho_r(h) &=& \left[ \max_{1\leq j \leq r} \left|\frac{ d^j}{dx^j}[w_{h, y}(x)]\right|_{x=y}\right]^{-1}.
\label{rho}
\end{eqnarray}

Let ${\cal G}_r$ be a class of functions $g(\theta)$ such that $p(y)$   defined by 
\fr{py}  is  $r$ times continuously differentiable in the neighborhood $\Omega_{y}$
of $y$. 
The following theorem provides lower bounds of the posterior risk \fr{rfun} at the point $y$.

\begin{theorem} \label{th:lower_bounds}
Let   $r \geq 1/2$,   Assumptions A1-A3 hold and $\rho_r(h)$ and $w_{h, y}(y)$ be such that,  for some 
nonnegative $r_1$ and $r_2$, 
\be \label{rhop}
\rho_r(h) \leq Ch^{r_1}, \quad |w_{h, y}(y)| \leq C_0h^{-r_2}.
\ee
Then, for any $y$ such that $p(y)>0$,  as $n \rightarrow \infty$,
\be \label{deltaupp1}
\Delta_n(y) = \inf_{\hat{t}} \  \sup_{g \in {\cal G}_r} \  \EE (\hat{t}  (y) - t(y))^2  
\geq  C n^{-\frac{2\max\{r, r_1\}-2r_2}{2\max\{r, r_1\} +1}},
\ee
where $C$ is an absolute constant independent of $n$.   In particular, if  $r_1 = r + r_2$ and $r_2= \max(\alpha/2,0)$
where $\alpha$ is defined in \fr{gamma_param_form}, then
\be \label{deltaupp}
\Delta_n(y)   \geq  C n^{-\frac{2r}{2r +1 + \max(\alpha, 0)}}  \quad   (n \to \infty).
\ee
\end{theorem}

\begin{remark}{\rm 
Note that if  $r_1 = r + r_2$ and $r_2= \max(\alpha/2,0)$, then EB estimator \fr{ebes} with $m=m_0$
given by formula \fr{m_opt} is asymptotically optimal. 
If $r_1 \leq r$ and $r_2=0$, then convergence rates are defined by behavior of  $p(x)$
in the neighborhood of $y$,  otherwise, the rates are defined by behavior of  $\Psi(x)$ in the neighborhood of $y$.
}
\end{remark}

\section{Adaptive choice of the resolution level using Lepski method} 
\label{sec:Lepskii_method} 
\setcounter{equation}{0}

 Note that  the value of $m_0$  depends on the unknown smoothness $r >0$ of  functions $p(x)$ and $\Psi(x)$ which is unknown,
i.e., the estimator is {\it nonadaptive}.  
In order to construct an adaptive estimator, we shall apply Lepski method
 for the optimal selection of the resolution level described  in  Lepski  (1991) and   
Lepski, Mammen and Spokoiny (1997).

Let $\gamma_m$ satisfy condition \fr{gamma_param_form}. 
Denote
\be \label{lambda_notation}
\rho^2_{mn}=   2^m (1 + \gamma^2_{m}) n^{-1} \log n,
\ee
and let $m_1$ and $m_n$ be such that
\be \label{eq:mn}
2^{m_1}= \log n,\ \ \ 
2^{m_n} ( \gamma^2_{m_n} +1) \asymp (\log n)^{-2} n.
\ee
Denote ${\cal M} = \{ m:\, m_1 \leq m \leq m_n\}$ and observe that 
under assumption \fr{gamma_param_form} with $b=0$,   $m_n$ is such that
$$
2^{m_n} = C\, \lkr \frac{n}{  \log^2 n} \rkr^{\frac{1}{1 + \max(\alpha,0)}}, 
$$
so that, for $m_0$ given by \fr{m_opt}, one has  $m_n/m_0 \rightarrow \infty$ as $n \rightarrow \infty$.

In this situation, the Lepski  method suggests to choose resolution level   as $m = \widehat{m}$ with
\begin{equation} \label{lepskii_eq}
\widehat{m}= \min \left\{ m \in {\cal M}: \left| \widehat{t}_{m}(y)-\widehat{t}_{j}(y)\right|^{2}
\leq \lambda^2 \left[ \|\hat{B}_{\delta m}^{-1} \|^2 + \|\hat{B}_{\delta j}^{-1} \|^2 \right]^2 
\rho_{jn}^2\ \ \mbox{for any}\ \ j\geq m \right\}, 
\end{equation}
where $\lambda$ is a positive constant independent of $m$  and $n$ and, for any $k$, 
$\hat{B}_{\delta k} = \hat{B}_{k} + \delta_k I$ where $\delta_k = 2^{k/2} n^{-1/2}$, 
and the matrix norm used in \fr{lepskii_eq} and throughout the paper is the spectral norm.

Let $\nu_1$ and $\nu_2$ be small positive values, $\nu_1 + \nu_2 <1$,  and 
$M$ be the size of vector $c$ and matrix $B$. Denote 
$C_{\varphi}= \displaystyle\sum\limits_{k} | \varphi(z-k)|$  and
\be \label{eq:S_value}
D= \|p\|_{\infty} + \| \Psi \|_{\infty}\| \varphi\|_{\infty} M\sqrt{M}[1-\nu_1]^{-1/2} + 
16\|\Psi\|_{\infty}\|\varphi\|^2_{\infty}\|p\|^{3/2}_{\infty}M^3 \sqrt{M}[1- \nu_2]^{-1}
\ee
and let $\lambda$ be given by equation
\be \label{eq:lambda_val}
\lambda = 16 C_{\varphi} \|p\|_{\infty}^{1/2} \sqrt{M} D  +1
\ee
Then, the following statement is true.

\begin{theorem}  \label{th:adaptive}
Let twice continuously differentiable scaling function $\varphi$ satisfy  \fr{wavco2} and \fr{wavco1}.
Let functions $p(x)$ and $\Psi(x)$ be $r \geq 1/2$   times continuously differentiable in the neighborhood 
$\Omega_{y}$ of $y$ and let $\Omega_{m,y} \subseteq \Omega_{y}$ where $\Omega_{m,y}$ is defined in \fr{ommy}.
Let $\gamma_m$ satisfy inequality \fr{gamma_param_form}. Construct EB estimator of the form \fr{ebes}
and choose $\widehat{m}$ according to \fr{lepskii_eq} with $\lambda$   
defined in \fr{eq:lambda_val} where $D$ given by \fr{eq:S_value}.  
If, for any $k \in K_{m,y}$,
\be \label{uniform_cond}
\| u_{m, k}\|_{\infty} \leq  C_u 2^{\frac{m}{2}}\gamma_m, 
\ee
then
\be \label{adapt_error}
 \EE \left| \widehat{t}_{\widehat{m}} (y) - t(y) \right|^2  
=  O \lkr n^{-\frac{2r}{2r+ 1 + \max(\alpha, 0)}}\, \log n \rkr.
\ee
\end{theorem}

In order to see how the method works, note that the   mean squared error can be decomposed  as
\be \label{delta}
\Delta = \EE \left| \widehat{t}_{\widehat{m}} (y) - t(y) \right|^2 = \Delta_1 + \Delta_2 
\ee
where $m_0$ is the optimal resolution level defined in formula \fr{mnopt}
and 
\begin{eqnarray*}
\Delta_1 & = & \EE[\left| \widehat{t}_{\widehat{m}} (y)- t(y) \right|^2 \II(\widehat{m}\leq m_0)],\\
\Delta_2 & = & \EE[\left|\widehat{t}_{\widehat{m}}(y)- t(y)  \right|^2 \II( \widehat{m} >m_{0} )].
\end{eqnarray*}
If $\widehat{m}\leq m_0$, then, by definition of $\widehat{m}$, one has
\be \label{eq:Delta1_portion}
\left|\widehat{t}_{\widehat{m}}(y)-\widehat{t}_{m_0}\right|^2\leq 
\lambda^2 \rho^2_{m_0 n}  \left[ \|\hat{B}_{\delta m}^{-1} \|^2 + \|\hat{B}_{\delta m_0}^{-1} \|^2 \right]^2 = O(\rho_{m_0 n}^2),
\ee
so that
\begin{eqnarray} \label{Delta1}
\Delta_1 & \leq & 2 \lkv \EE\left|\widehat{t}_{\widehat{m}}(y)-\widehat{t}_{m_0}(y) \right|^2 
+ \EE \left| \widehat{t}_{m_0}(y) - t(y) \right|^2 \rkv = O \lkr \rho^2_{m_0n}  \rkr.
\end{eqnarray}
Now, in the case when $\h{m}> m_0$,  by definition \fr{lepskii_eq} of $\widehat{m}$, there  exists   
$j  > m_0$,   such that
$$
| \h{t}_j(y)-\h{t}_{m_0}(y) |^2 > \lambda^2 \rho^2_{j n}  \left[ \|\hat{B}_{\delta j}^{-1} \|^2 + \|\hat{B}_{\delta m_0}^{-1} \|^2 \right]^2 . 
$$ 
It turns out that probability of such an event is very low. In particularly, the following lemma holds.

\begin{lemma} \label{lem:prob4} 
Let conditions of Theorem \ref{th:adaptive} hold. 
If resolution level $m$  is such that $m_0 < m \leq m_n$, then, as $n \rightarrow \infty$,   
\be \label{prob4}
\PP \lkr | \h{t}_m(y)-\h{t}_{m_0}(y) |^2 > \lambda^2 \left[ \|\hat{B}_{\delta m}^{-1} \|^2 + 
\|\hat{B}_{\delta m_0}^{-1} \|^2 \right]^2  \rho^2_{mn} \rkr =O\left( n^{-2}\right).
\ee
where $\rho^2_{mn}=   2^m (1 + \gamma^2_{m}) n^{-1} \log n$.
\end{lemma}

This lemma implies that $\Delta_2 = O(n^{-1}) = o\lkr \rho^2_{m_0n}  \rkr$ as $n \rightarrow \infty$
and, hence, Theorem \ref{th:adaptive} is valid. The complete proof of  Theorem \ref{th:adaptive} 
is provided in Section \ref{sec:proofs}.

\section{Examples of construction of EB estimators, lower and upper bounds for the risk} 
\label{sec:examples} 
\setcounter{equation}{0}

\subsection{Location Parameter Family}

In the case when $\theta $ is a location parameter, one has  $q(x|\theta)  =  q(x- \theta)$,  
EB estimator  $t(y)$ in \fr{pm}  is of the form
$$
t(y)= y - \frac{\int^{\infty}_{-\infty} (y- \theta) q(y- \theta)g(\theta) d \theta}{\int^{\infty}_{-\infty} q(y- \theta) g(\theta) d \theta}.
$$
Hence, in this case, $u_{m, k}(x)$  is the solution to the following equation
\be \label{eq Loc1}
\int^{\infty}_{-\infty} q(x-\theta) u_{m, k}(x)dx= \int^{\infty}_{-\infty} (x- \theta)q(x-\theta) \varphi_{m, k}(x)dx.
\ee
Here we provide a brief construction  of $u_{m, k}$  and evaluation of its variance, for a more detailed derivation 
see Pensky (2000, 2002). Taking Fourier transforms of both sides, we obtain that 
\be  \label{eq:Umk}
u_{m, k}(x)= 2^{m/2}U_m(2^mx-k),
\ee
where $U_m (x)$   is the inverse Fourier transform of 
\be  \label{eq:Umomega}
\h{U}_{m}( \omega) =   i  [\h{q}(- 2^{m} \omega)]^{-1}\,    \h{ q}'(- 2^{m} \omega)\, \h{\varphi}(\omega). 
\ee
 Here, and in what follows,  $\hat{f}$  denotes the Fourier transform of   function $f$.  
Using  Parseval identity and taking into account that $K_{m,y}$ has finite number of terms, we derive  
\be \label{Parseval}
\gamma_m^2 \asymp \int^{\infty}_{-\infty}  | \h{q}(- 2^m \omega)]^{-1}   
\h{ q}'(-2^m \omega) |^2  | \h{\varphi}( \omega)|^2 d \omega.
\ee
Also, the following relation allows to check validity of condition  \fr{uniform_cond} 
\be \label{loc_uniform}
\|u_{m, k}(x)\|_{\infty} \leq 2^{m/2} \|U_m(x)\|_{\infty}.
\ee

Now, in order to calculate the lower bounds for the risk, 
we need to find $\psi_{h, y}(\theta)$ and $w_{h, y}(x)$. 
Let $\psi_{h, y}(\theta)$ be solution of equation \fr{psi-equation}.
It is easy to show that $\psi_{h, y}(\theta)$ is of the form 
$ \psi_{h, y}(\theta)= \psi_h ((\theta -y)/h),$
where the Fourier transform $\h{\psi}_h(\omega)$ of $\psi_h( . )$ is of the form
\be  \label{eqpsif}
\h{\psi}_h(\omega) =  \h{k}(\omega)/\h{q}(\omega/h).
\ee
 In order to obtain an expression for $w_{h, y}(x)$, recall that equation \fr{psi-equation}  can be rewritten as
\begin{eqnarray*}
w_{h, y}(x) & =& x- \int_{-\infty}^\infty  \left( x-\theta \right)q(x- \theta) \psi_{h, y}(\theta)d\theta 
= x-w_h\left(\frac{x-y}{h}\right),
\end{eqnarray*}
where  $w_h( . )$ is the inverse Fourier transform of 
\be  \label{eqomf}
\h{\omega}_h(\omega)= i^{-1}\h{k}(\omega) {\h{q}'(\omega/h)}/{\h{q}(\omega/h)}.
\ee
In this situation,   $\rho_r(h)$ defined in \fr{rho} and $w_{h, y}(y)$ are, respectively,  given by
$$
w_{h, y}(y)  = y-w_h(0), \quad
\rho_r(h) =  \left[ \max_{1\leq j \leq r} \left(h^{-j}[\omega^{(j)}_{h}(0)]\right)\right]^{-1} 
$$
Below, we  consider some special cases.
\\

\begin{example} \label{Nordi}
{\bf Normal Distribution  }  
{\em Let $q(x|\theta)$ be the pdf of the normal distribution  
$$
q(x|\theta) = \frac{1}{ \sqrt{2 \pi} \sigma} e^{-\frac{ (x-\theta)^2}{2 \sigma^2}},
$$  
where $\sigma >0$ is known. Then 
$\hat{q}( \omega)=\sqrt{ \pi} e^{-\frac{ \omega^2 \sigma^2}{2}}$ 
and, using properties of Fourier transform, we derive
$$
U_m(x)= - 2^m \sigma^2 \varphi'(x).
$$
Since $\varphi'(x)$ is square integrable,   it follows from \fr{Parseval} and \fr{loc_uniform} that
 $ \gamma^2_m \asymp 2^{2m}$  and $\|U_m\|_{\infty} \asymp 2^m$. Hence, condition \fr{uniform_cond} of Theorem
\ref{th:adaptive} holds. Then,   $\alpha =2$,  $2^{m_0}\sim n^{\frac{1}{2r+3}}$, and   
$$
 \EE \left| \widehat{t}_{\widehat{m}} (y) - t(y) \right|^2  
=  O \lkr n^{-\frac{2r}{2r+ 3}}\, \log n \rkr.
$$
by Theorem \ref{th:adaptive}.  

In order to verify optimality of the estimator, we derive the lower bounds for the risk following Theorem~\ref{th:lower_bounds}.
Using \fr{eqomf}, by direct calculations, we  obtain, $\h{\omega}_h(\omega)=  i h^{-1} \omega \sigma^2\, \h{k}(\omega)$.
Consequently, for a fixed value of $y$, one has 
$w_{h, y}(x) = x+\sigma^2h^{-1} k'(h^{-1} (x-y))$, 
$\rho_r(h) =  C h^{r+1}$ and $|w_{h, y}(y)| \leq  C h^{-1}$.
Hence, $r_1= r+r_2$, $r_2= \alpha/2 = 1$  in \fr{rhop}, and 
   application of  Theorem~\ref{th:lower_bounds} yields 
$\Delta_n(y)   \geq  C n^{-\frac{2r}{2r+ 3}}$, so that EB estimator is optimal up to a logarithmic factor.
}
\end{example}

\begin{remark}  \label{remark2}
{\rm 
Note that in the case of the normal distribution, 
\begin{eqnarray*}
p^{(r)} (x) & = & \int_{-\infty}^{\infty} \frac{1}{ \sqrt{2 \pi} \sigma}  
{\cal P}_r \lkr \frac{x-\theta}{\sigma} \rkr\,   e^{-\frac{ (x-\theta)^2}{2 \sigma^2}} \ g(\theta) d \theta, \\
\Psi^{(r)} (x) & = & \int_{-\infty}^{\infty} \frac{1}{ \sqrt{2 \pi} \sigma}  
{\cal P}_r \lkr \frac{x-\theta}{\sigma} \rkr\,   e^{-\frac{ (x-\theta)^2}{2 \sigma^2}} \ \theta \, g(\theta) d \theta,  
\end{eqnarray*} 
where ${\cal P}_r (z)$ is the polynomial in $z$ of degree $r$. Since ${\cal P}_r (z) e^{-z^2/2}$ is an infinitely smooth uniformly bounded function, 
$r$ can take limitlessly large values, as long as $|\theta| \, g(\theta)$ is integrable. 
Then, convergence rate of the EB estimator is  bounded  only by the number 
of vanishing moments of the scaling function $\varphi$. 
}
\end{remark}

\begin{example} \label{dblexp}
{\bf Double-exponential distribution} 
{\em  
   Let $q(x|\theta)$ be the pdf of the double-exponential distribution
$$
q(x|\theta) =\frac{1}{2 \sigma} e^{ -\frac{|x-\theta|}{\sigma}},
$$  
where $\sigma >0$ is known. Then
$\hat{q}( \omega)=  1/(1+  \omega^2 \sigma^2)$, and, using properties of Fourier transform, we obtain 
$$
u_{m, k}(x)= \int ^{\infty}_{- \infty} \varphi_{m, k}(t) \mbox{sign} (x-t) \exp( -|x-t|/\sigma) dt.
$$
Note that for any $\omega$, one has $\left| \hat{q}(-  \omega)]^{-1}   \hat{ q}'(- \omega) \right|  \leq \sigma$. Therefore, 
$$ 
\gamma^2_m \asymp \int^{\infty}_{-\infty}  \left| \hat{\varphi}( \omega) \right|^2 d \omega=  1,\ \ \ \ 
|U_m(x)| \leq \int^{\infty}_{-\infty} \left| \hat{\varphi}(\omega) \right|d\omega  \asymp 1,
$$
so that condition \fr{uniform_cond} holds. Consequently, $\alpha=0$, $2^{m_0}\sim n^{\frac{1}{2r+1}}$, and 
$$
 \EE \left| \widehat{t}_{\widehat{m}} (y) - t(y) \right|^2  
=  O \lkr n^{-\frac{2r}{2r+ 1}}\, \log n \rkr.
$$
by Theorem \ref{th:adaptive}.

In order to derive  lower bounds for the risk,  we apply \fr{eqomf} to obtain  
$$
\h{\omega}_h(\omega)= - i^{-1} \h{k}(\omega) \ \frac{2\omega \sigma^2h^{-1}}{\sigma^2\omega^2h^{-2} +1}. 
$$
It can be shown that 
\begin{eqnarray*} 
w_{h, y}(x)&=&x-2   \int^1_{-1}  k\left( t\right) \mbox{sign} (x-y-ht) \exp( - \sigma^{-1}\, |x-y-ht|) dt \\
\end{eqnarray*}
$\rho_r(h) = C h^{-1}$  and
$$
| w_{h, y}(y)| = y+ 2h \left| \int^1_{-1} k(t) e^{- \frac{h |t|}{\sigma}} dt \right|.
$$
Therefore,  $r_1= r$ and $r_2=0$ in \fr{rhop}, and    application of  Theorem~\ref{th:lower_bounds} yields 
$\Delta_n(y)   \geq  C n^{-\frac{2r}{2r+ 1}}$, so that EB estimator is optimal up to a logarithmic factor.
} 
\end{example}

\subsection{One-Parameter exponential Family}

Let conditional  distribution belong to a one-parameter exponential family, i.e. 
$$
q(x\mid \theta)= h( \theta) f(x) e^{-{x^{b }}{ \theta}},  \ x \in {\cal X}, \ \theta \in \Theta, 
$$
where $h( \theta) > 0$ and $f(x) > 0$. Then, equation \fr{eq1} is of the form 
\begin{eqnarray*}
\int_{\cal X} f(x)  e^{-{x^{b }}{ \theta}}u_{m, k}(x)dx &=& 
\int_{\cal X} \theta f(x) e^{-{x^{b }}{ \theta}}\varphi_{m, k}(x)dx. 
\end{eqnarray*}
 Integrating the right hand side  by parts and solving for $u_{m, k}$ we derive 
 \be \label{uscale1}
u_{m, k}(x)= \frac{1}{b  f(x)} \ \frac{d}{dx}\lfi \frac{ f(x)\varphi_{m, k}(x)}{x^{b -1}} \rfi 
=  \frac{2^{3m/2} \varphi'(2^m x-k)}{b  x^{b -1}} + \frac{(f(x) x^{1 - b })'  \, 2^{m/2} \varphi (2^m x-k)}
{b   f(x)}.
\ee
Let the value of $y$ be such that $c_1 \leq y \leq c_2$ for some $0< c_1 < c_2 < \infty$.
Then, it is easy to show that, if $k \in K_{m,y}$, then $k \asymp 2^m$. In this case, 
for any $k \neq 0$, one has
\begin{eqnarray}  \label{gamma-component}
[\gamma^{(1)}_k(m)]^2 & \leq & 2 b ^{-2}\, 2^{2b  m} \   \int_{M_1}^{M_2} |z+k|^{-(2b  -2)} [\varphi^\prime (z)]^2  dz \\
 & + & 2 b ^{-1} \max_{M_1 \leq 2^m x - k \leq M_2} \left|(f(x))^{-1} (f(x) x^{1 - b })'\right| \leq C 2^{2m}, \nonumber
\end{eqnarray} 
so that $\alpha =2$.

 In order to evaluate lower bounds  for the risk,  we need to find $\psi_{h, y}(\theta)$ and $w_{h, y}(x)$. 
Let $\psi_{h, y}(\theta)$ be solution of equation \fr{psi-equation} and $w_{h, y}(x)$ be defined by \fr{w}.
 It is straightforward to verify that 
 \begin{eqnarray}  \label{us}
   w_{h, y}(x) & = & - \frac{f(x)}{b  x^{b  - 1}}\ \frac{d}{dx} \left[ \frac{1}{f(x)}\, k \lkr  \frac{x-y}{h}\rkr \right] \\
 & =   &
\frac{f'(x)}{b  x^{b  -1} f(x)}   k \lkr  \frac{x-y}{h}\rkr - \frac{1}{b  h x^{b  -1}}   k' \lkr  \frac{x-y}{h}\rkr.
\nonumber 
\end{eqnarray} 
For a fixed value of $y>0$, one has  $r_1= r +1$ and $r_2=1$ in \fr{rhop}. Hence, application of  Theorem~\ref{th:lower_bounds} 
yields  
\be \label{onepexp}
\Delta_n(y)   \geq  C n^{-\frac{2r}{2r+ 3}}.
\ee

 \begin{example} \label{Weib} 
 {\bf Weibull Distribution   }  
{\em 
   Let $q(x| \theta)$ be the pdf of the  Weibull distribution  
 \begin{eqnarray*}
q(x| \theta)= {b }{\theta}x^{b -1}e^{-{x^{b }}{\theta}}, \ x \geq 0, \ \theta >0,\ b  \geq 1.
 \end{eqnarray*}
 In this case,   $f(x)= x^{b -1}$, $h(\theta)= {b }{\theta}$  and  $u_{m, k}(x)$ is of the form \fr{uscale1}
with the second term being identical zero 
\begin{eqnarray} \label{Weiu}
u_{m, k}(x) =  \frac{2^{3m/2} \varphi'(2^mx-k)}{b  x^{b -1}}.   
\end{eqnarray}
Hence, it follows from formula \fr{gamma-component} that 
\begin{eqnarray}
[\gamma^{(1)}_k(m)]^2  \leq 
2^{2b  m} b ^{-2} (M_2 - M_1)  |M_1 + k|^{-(2b  -2)}.
\label{Wei:gammajm}
\end{eqnarray}
Since $k \asymp 2^m$ and  the set $K_{m,y}$ has a finite number of terms,
\fr{Wei:gammajm} yields that $\gamma^2 _m \asymp 2^{2m}$.

 Finally, it remains to verify whether the condition \fr{uniform_cond} of Theorem \ref{th:adaptive} holds.  
Indeed, it follows from  \fr{Weiu} that for $k \in K_{m,y}$, one has 
\be
\sup_x| u_{m, k}(x)| \leq \sup_z \left| \frac{ 2^{3m/2} \varphi'(z) 2^{m(b -1)}}{b  [ M_1+ k]^{b  -1}} \right| 
\leq C 2^{m/2} \gamma_m, 
\ee
so Theorem  \ref{th:adaptive} can be applied.
Hence,   under   assumptions of  Lemma\ \ref{lem:syserror}, 
one has $\alpha  = 2$, $2^{m_0}\sim n^{\frac{1}{2r+3}}$, and, by Theorem  \ref{th:adaptive},
$ 
\EE \left| \widehat{t}_{\widehat{m}} (y) - t(y) \right|^2  
=  O \lkr n^{-\frac{2r}{2r+ 3}}\, \log n \rkr.
$ 
Therefore,  the EB estimator is optimal within a log-factor of $n$ due to \fr{onepexp}. 
}
  \end{example}

 \begin{example}\label{Gamdis}
 {\bf Gamma Distribution. }  
{\em
Let the pdf of the  Gamma distribution be given by
  \begin{eqnarray} \label{gammapdf}
q(x| \theta)= {\frac{1}{\Gamma{(\beta)}}}{\theta^{\beta}}x^{\beta-1}e^{-{x}{\theta}}, \ x \geq 0, \ \theta >0,\ \beta \geq 1.
 \end{eqnarray}
In this case, $b =1$, $f(x) = x^{\beta -1}$ and $h(\theta) = \theta^{\beta}/\Gamma(\beta)$.
Note that the family of densities  \fr{gammapdf} is also a scale parameter family.
By formula \fr{uscale1},  $u_{m, k}(x)$ is of the form
\begin{eqnarray*}
u_{m, k}(x)&=&  (\beta -1 ) x^{-1} 2^{m/2} \varphi(2^mx-k)  + 2^{3m/2} \varphi'(2^mx-k).
 \end{eqnarray*}
If $y$ is such that $c_1 \leq y \leq c_2$ for some $0< c_1 < c_2 < \infty$, then, by calculations similar to Example~\ref{Weib},
obtain that   $2^{m_0}\sim n^{\frac{1}{2r+3}}$ and  Theorem \ref{th:adaptive} holds, so that,
$ 
\EE \left| \widehat{t}_{\widehat{m}} (y) - t(y) \right|^2  
=  O \lkr n^{-\frac{2r}{2r+ 3}}\, \log n \rkr.
$ 
Therefore,  the EB estimator is optimal within a log-factor of $n$ due to \fr{onepexp}. 
}
  \end{example}

  \subsection{Scale parameter family} 

If $q(x|\theta)$ is a scale parameter family, 
$q(x|\theta) =  \theta^{-1}\, q \lkr \theta^{-1} x  \rkr$, 
it is difficult to pinpoint a general rule for finding $u_{m,k} (x)$,
however, as it follows from Example~\ref{Gamdis}, many particular cases 
can be treated. Below, we consider   one more example.

 \begin{example}\label{Uniform}
 {\bf Uniform Distribution  }  
{\em
Let $q(x \mid \theta)$ be given by
  $$
q(x \mid \theta) = \theta^{-1}\, \II \left( 0 < x< \theta\right), \  \  a \leq \theta \leq b. 
$$
  Then, equation \fr{eq1} is of the form 
   \begin{eqnarray} \label{unif1}
   \int ^{\theta}_0  \theta^{-1}  u_{m, k}(x) dx = \int^{\theta}_0 \varphi_{m, k} (x) dx
    \end{eqnarray}
Taking derivatives  of both sides of \fr{unif1}  with respect to $\theta$ and replacing $\theta$ by $x$, we derive
       \begin{eqnarray}
   u_{m, k} (x)= 2^{-m/2} \int^{2^mx-k}_{M_1} \varphi(z)dz + x 2^{m/2} \varphi(2^mx-k).    
           \end{eqnarray}
Since $a \leq \theta \leq b$, then also $a \leq x \leq b$, and it is easy to check that
\begin{eqnarray*}
\int_a^b  x^2 2^{m} \varphi^2 (2^mx-k) dx \asymp 1,\ \ \ 
\int_a^b  \lkr 2^{-m/2} \int^{2^mx-k}_{M_1} \varphi(z)dz \rkr^2 dx = O(2^{-m}),
\end{eqnarray*}
as $m \rightarrow \infty$. Then,  $\gamma_m \asymp 1$, $\alpha= 0$ and condition \fr{uniform_cond} holds.
Therefore,   $2^{m_0}\sim n^{\frac{1}{2r+1}}$  and, by Theorem  \ref{th:adaptive},  one has
$\EE \left| \widehat{t}_{\widehat{m}} (y) - t(y) \right|^2  
=  O \lkr n^{-\frac{2r}{2r+ 1}}\, \log n \rkr.$
\\

 Now, in order to calculate the lower bound for the risk,  we need to find $\psi_{h, y}(\theta)$ and $w_{h, y}(x)$. 
According \fr{psi-equation}  and \fr{w},  functions $\psi_{h, y}(\theta)$ and  
$w_{h, y}(x)$ satisfy equations 
 \begin{eqnarray*}
\int^{b}_{x}  \theta^{-1} \psi_{h, y}(\theta)d\theta  & = &  k\left(h^{-1} (x-y) \right), \quad 
\int^{b}_{x}  \psi_{h, y}(\theta)d\theta =  w_{h, y}(x).
\end{eqnarray*}
Now, taking derivatives of  both sides of the first equation with respect to $x$ and solving for 
$\psi_{h, y}(\theta)$, we obtain
$$
\psi_{h, y}(\theta)=- h^{-1}  \theta k'\left(h^{-1} (x-y) \right).
$$
It can be shown that 
$$
w_{h, y}(x)= x k(h^{-1} (x-y)) + h K(h^{-1} (x-y)),
$$
where  $K'(z)=k(z)$. Note that $r_1= r$ and $r_2=0$, hence,  applying Theorem \ref{th:lower_bounds},
we obtain the following lower bounds for the risk   
$\Delta_n(y)   \geq  C n^{-\frac{2r}{2r+ 1}}$, so that the EB estimator is optimal, up to a logarithmic factor.
}
\end{example}

\section{Discussion}
\label{sec:discussion}
\setcounter{equation}{0}

The present paper achieves  two main objectives. The first one is to derive lower bounds for the postrior risk
of a nonparametric empirical Bayes estimator  under general assumptions. The present paper is the first one to accomplish this task.
The second purpose of this paper is to provide an adaptive  wavelet-based method of  EB estimation.
The method is based on approximating   Bayes estimator $t(y)$ corresponding to observation $y$ 
as a whole using finitely supported wavelet family. The wavelet estimator is used 
in a rather non-ortodox way: $t(y)$ is estimated locally using 
only a linear scaling part of the expansion at the resolution level $m$ where 
coefficients are recovered by solving a system of linear equations.

The advantage of the method lies in its flexibility.  The technique works for a variery of families 
of conditional distributions. Computationally, it leads to solution of a
finite system of linear equation which, due to decorrelation property of wavelets, is sparse and well conditioned. 
The size of the system depends on a size and regularity of the wavelet which is used for representation
of the EB estimator $t(y)$.

A non-adaptive version of the method was introduced in Pensky and Alotaibi (2005).
However, since no  mechanism for choosing the resolution level $m$ of the expansion 
was suggested, the Pensky and Alotaibi (2005) paper remained of a theoretical interest  only.
In the present paper, we use Lepski method  for choosing an optimal resolution level 
$m$ and show that the resulting EB estimator remains nearly asymptotically optimal (within a 
logarithmic factor of the number of observations $n$).  We also show that the EB estimators 
constructed in the paper are asymptotically optimal (Up to a logarithmic factor) as $n \to \infty$.

Finally, we should comment that, although the choice of a wavelet basis for representation of 
$t(y)$ is convenient, it is not unique. Indeed, one can use a local polynomial  or a kernel 
estimator  for representation of $t(y)$. In this case, the challenge of finding support of the estimator 
for the local polynomials or bandwidth for a kernel estimator can be addressed by Lepski
method in a similar  manner. However, the disadvantage of abandoning wavelets will be that 
the system of equations will cease to be sparse and well-posed.

\section{Proofs} 
\label{sec:proofs}
\setcounter{equation}{0}

 In what follows, we suppress index $m$ in notations of matrix $B_m =B$, $B_{m\delta} =B_{\delta}$, $\h{B}_m =\h{B}$ and $\h{B}_{m\delta} =\h{B}_{\delta}$, and vector $c_m =c$
unless this leads to a confusion. 
Proofs of Lemma  \ref{lem:randerror} is based on  the following lemmas.

\begin{lemma}  \label{lem:hBhc}
Let  $B$, $c$, $\h{B}$ and $\h{c}$   be defined  in  \fr{bjk},  \fr{hbjk}, \fr{cj} and \fr{hcj}, respectively, and let $M$ be the size of the vector $c$.
If  $n \rightarrow \infty$, $2^m/n \rightarrow 0$, one has 
\be \label{hBhc1}
\EE \left|\left|\h{B}-B\right|\right|^{2l} = O\left( n^{-l} 2^{ml}\right), \ \  l =1, 2, 4,  
\ee 
\be  \label{hBhc2}
\EE \left|\left|\h{c}-c\right|\right|^{2} = O\left(n^{-1} \gamma^2_m \right),\\
\ee
\be  \label{hBhc3}
\EE \left|\left|\h{c}-c\right|\right|^{4} = O\left(n^{-3}\left|\left|\gamma^{(2)}(m)\right|\right|^2 + n^{-2} \gamma^4_m \right),
\ee
and
\begin{eqnarray} 
\EE \left|\left|\h{c}-c\right|\right|^{8} &=& O\left({n^{-7}}{\left|\left|\gamma^{(4)}(m)\right|\right|^2}
+{n^{-6}}{\left|\left|\gamma^{(2)}(m)\right|\right|^4}+ {n^{-6}}{\gamma^2_m\left|\left|\gamma^{(3)}(m)\right|\right|^2}\right)
\nonumber \\
&+&O\left({n^{-5}} {\left|\left|\gamma^{(1)}(m)\right|\right|^4\left|\left|\gamma^{(2)}(m)\right|\right|^2}
+{n^{-4}}{\left|\left|\gamma^{(1)}(m)\right|\right|^8}\right).  \label{hBhc4} 
\end{eqnarray} 
If, in addition, \fr{uniform_cond} holds, then, as $n \rightarrow \infty$, 
\be \label{c_dev_moments}
\EE  \| \h{c}-c \|^{2\kappa} = O \lkr n^{-\kappa}  2^{m(\kappa-1)} \gamma_m^{2\kappa} \rkr, \kappa=1,2,4.
\ee
\end{lemma}

{\bf Proof of Lemma  \ref{lem:hBhc}.} 
Recall that $\hat{B}_{j,k}-B_{j,k}= n^{-1} \displaystyle\sum_{t=1}^{n} \eta_t$
where $\eta_t = \varphi_{m,k} (X_t) \varphi_{m,j} (X_t)$, $t=1, \cdots, n$. Hence, taking the second moment we obtain 
\begin{eqnarray*}
\EE \left|\left|\h{B}-B\right|\right|^2& \leq& M^2  n^{-1} \EE [\eta^2_t]  
 \leq n^{-1}  M^2 \left[ 2 \left|\left|p\right|\right|_{\infty}\left|\left|\varphi\right|\right|^2_{\infty}2^m \right]  
 = O \left( n^{-1}  2^m \right).
\end{eqnarray*}
For $l=2$, we  apply Jensen's inequality 
\begin{eqnarray*}
\EE\left[ n^{-1} \displaystyle\sum_{t=1}^{n} \eta_t\right]^4&=& O \lkr n^{-4}  \lkv n\EE [\eta^4_1]+ n(n-1)\EE^2 [\eta^2_1] \rkv \rkr  
= O \lkr n^{-2} \EE [\eta^4_1] \rkr.
\end{eqnarray*}
Since $\EE [\eta^4_1] \leq  2^{2m} \| p \|_\infty \|\varphi \|^2_{\infty}$ and matrix $B$ is of  finite dimension $M$,
\fr{hBhc1} is valid for $l=2$. In a similar manner we can show that \fr{hBhc1}  holds for $l=4$.

In order  to prove \fr{hBhc2}--\fr{hBhc4}, recall that 
$\h{c}_k-c_k= n^{-1}\, \displaystyle\sum^n_{t=1} \xi_t$
where $\xi_t = u_{mj}(X_t)$, $t=1, \cdots, n$. Thus,  
\begin{eqnarray*}
\EE \|\h{c}-c\|^2 & \leq & M  n^{-1} \EE [\xi^2_1]  
 \leq  M  n^{-1}  \|p \|_{\infty}\gamma^2_m   
=  O\left( n^{-1} \gamma^2_m  \right)
\end{eqnarray*}
which proves \fr{hBhc2}. 
Now, to prove \fr{hBhc3}, observe that 
\begin{eqnarray*}
\EE\left[ n^{-1}\, \displaystyle\sum_{t=1}^{n} \xi_t\right]^4 & = & O \lkr n^{-3} \EE [\xi^4_1]+ n^{-2} \EE^2 [\xi^2_1] \rkr 
=  O\left(n^{-3}\left|\left|\gamma^{(2)}(m)\right|\right|^2+n^{-2}\left|\left|\gamma^{(1)}(m)\right|\right|^4\right) 
\end{eqnarray*}
and note that vector $c$ is of   finite dimension $M$. 
The proof of \fr{hBhc4} can be carried out in a similar manner.

Now, let us check validity of \fr{c_dev_moments} when condition \fr{uniform_cond} holds. 
Observe that, under condition \fr{uniform_cond},  for   $\kappa=1, 2,  ....$ one has
\be  \label{jnormbnd}
\| \gamma^{(\kappa)}(m) \|^2  \leq   C^{2\kappa-2}_u 2^{m(\kappa-1)} \gamma^{2\kappa}_m 
\ee
Plugging \fr{jnormbnd} into \fr{hBhc2}--\fr{hBhc4}, obtain \fr{c_dev_moments}.

\begin{lemma}  \label{lem:large_dev} 
Let  $B$, $c$, $\h{B}$ and $\h{c}$   be defined  in  \fr{bjk},  \fr{hbjk}, \fr{cj} and \fr{hcj}, respectively,
and let $M$ be the size of the vector $c$. If $2^m \leq n (\log n)^{-2}$, then, for any $\tau >0$,
\be \label{large_devB}
\PP\left(\left|\left|\h{B}-B\right|\right|^2\geq M^2 \tau^2 2^m \,n^{-1} \log {n}  \right)\leq 
2 M^2 n^{- \frac{\tau^2}{8\left|\left|\varphi\right|\right|^2_\infty\left|\left|p\right|\right|_\infty}}. 
\ee
If, in addition, \fr{uniform_cond} holds, then, for any $\tau >0$, 
\be \label{large_devc}
\PP\left(\left|\left|\h c-c\right|\right|^2\geq M \tau^2 \gamma^2_{m}  n^{-1} \log n \right)\leq 
2M n^{-\frac{\tau^2}{8 \left|\left|p\right|\right|_{\infty}}},
\ee
\end{lemma}

{\bf Proof of Lemma \ref{lem:large_dev} } 
The proof is based on application of Bernstein inequality
$$ P \lkr \left| n^{-1} \sum_{t=1}^n Y_t \right| > z \rkr \leq 
2\ \exp \lkr-\frac{nz^2}{2(\sigma^2 + \|Y\|_\infty z/3)}\rkr, $$
where  $Y_t,\  t=1,...,n,$ are i.i.d. with ${\bf E} Y_t =0$, 
 ${\bf E} Y^2_t = \sigma^2$ and $\|Y_t\| \leq \|Y\|_\infty < \infty$.

Recall that $\h{B}_{j,k}$ defined by \fr{hbjk} are the unbiased estimators of $B_{j,k}$
defined in \fr{bjk}. Denote 
$\eta_t=\varphi_{mj}(X_t)\varphi_{mk}(X_t)-B_{j,k}$, so that 
$\hat{B}_{j,k}-B_{j,k}=n^{-1}\, \displaystyle\sum_{t=1}^{n} \eta_t$, where 
$\eta_t$ are i.i.d. with $\EE\left(\eta_{t}\right)= 0$  and 
$\sigma^2 = \EE\left(\eta^2_t\right)\leq  
2\left|\left|p\right|\right|_\infty\left|\left|\varphi\right|\right|^2_\infty2^m$.
Also, 
$\left|\left|\eta_t\right|\right|_\infty \leq  2^{m+1} \left|\left|\varphi\right|\right|^2_\infty$. 

Applying Bernstein inequality with $z= \tau 2^{\frac{m}{2}}\sqrt{\log n /n }$
for every $k \in K_{m,y}$ and recalling that $2^m \log n/n \rightarrow 0$ as $n \rightarrow \infty$,  we obtain 
\begin{eqnarray*}
\PP \left(\left|\h B_{j,k}-B_{j,k}\right| >  \tau2^{\frac{m}{2}} \sqrt{\log n/n} \right)
& \leq &
2\exp\left(\frac{-\tau^2\log n}{4\left|\left|\varphi\right|\right|^2_\infty\left(\left|\left|p\right|\right|_\infty+2^{\frac{m}{2}}
\frac{\tau \sqrt{\log n}}{3\sqrt{n}}\right)}\right)\\
& \leq &  
2 \exp \left( \frac{-\tau^2\log n}{ 8\left|\left|\varphi\right|\right|^2_\infty \left|\left|p\right|\right|_\infty} \right) 
\end{eqnarray*}
Since matrix $B$ has $M^2$ components, \fr{large_devB} is valid.\\

Now, in order  to prove \fr{large_devc}, recall that $\h{c_j}$
given by \fr{hcj} is an unbiased estimator of $c_j$,
so that, for any $k$, variables $\xi_t = u_{m,k}(X_t) - c_k$, $t=1, \cdots, n$, are i.i.d. 
with $\EE \xi_t=0$ and 
$$
\EE\left(\xi^2_t\right)\leq  \int^{\infty}_{-\infty} 
u^2_{m,k}(x)p(x)\,dx \leq \left|\left|p\right|\right|_\infty \gamma^2_m.
$$ 
In addition, $\left|\left|\xi_t\right|\right|_\infty \leq2\left|\left|u_{m,k}\right|\right|_\infty  \leq 2 \gamma_m 2^{\frac{m}{2}}.$  
Thus, applying  Bernstein inequality with $z= n^{-1/2} \tau \gamma_{m} \sqrt{\log n}$, we obtain  
\begin{eqnarray*}
\PP\left(\left|\h{c}_k- c_k\right|> n^{-1/2} \tau \gamma_m \sqrt{\log n} \right)
& \leq &
2 \exp \left( \frac{ -\tau^2\log{n}}{2(  \left|\left|p\right|\right|_{\infty} +2\frac{ 2^{\frac{m}{2}} \tau \sqrt{\log n}}{3\sqrt{n}})}\right)
\leq
2 \exp \left( \frac{ -\tau^2\log{n}}{8 \left|\left|p\right|\right|_{\infty}}\right).
\end{eqnarray*}
To  complete  the proof, note that vector $c$ has $M$ components.\\




\begin{lemma}  \label{lem:adelta_dev}
Let $\Omega_B = \left\{ \omega: \|\hat{B} - B\| > 0.5\, \|B^{-1}\|^{-1} \right\}.$ Then,  
\beqn   \label{eq:hata_a}
\|  \hat{a}_{\delta} - a \| & \leq & \| \hat{c}- c\| \|B^{-1}\|   + \|c\| \lkv \delta \|B^{-1}\|^2 + 
2 \|B^{-1}\|^2 \| \hat{B}-B\| + 2 \delta^{-1}  \II(\Omega_B) \rkv \nonumber \\
& + & \| \hat{c}- c\| \lkv \delta \|B^{-1}\|^2 + 
2 \|B^{-1}\|^2 \| \hat{B}-B\| + 2 \delta^{-1}   \II(\Omega_B) \rkv.
\eeqn
\end{lemma}

{\bf Proof of Lemma  \ref{lem:adelta_dev}. }  Since $a=B^{-1}c$ and $\h{a}_{\delta } = (\h{B} + \delta I)^{-1} \h{c}$, 
by the properties of the norm, we obtain
\beqn
  \|  \hat{a}_{\delta } - a \| & \leq &  \| B^{-1}\| \| \hat{c}- c\| + \|\hat{B}^{-1}_{\delta } - 
B^{-1} \|\|c\| + \| \hat{B}^{-1}_{\delta } - B^{-1}\|\| \hat{c} - c \|
  \label{hadeltan}\\
 \| \hat{B}^{-1}_{\delta } - B^{-1} \| & \leq&  \| \hat{B}^{-1}_{\delta }- B^{-1}_{\delta } \| + \| B^{-1}_{\delta } - B^{-1} \|
 \label{hBdelta}
\eeqn
Now, since $\hat{B}^{-1}_{\delta }- B^{-1}_{\delta } = \hat{B}^{-1}_{\delta } ( \hat{B}-B)  B^{-1}_{\delta }$, 
$\| B^{-1}_{\delta }\| \leq \delta^{-1}$, $\|  \hat{B}^{-1}_{\delta }\| \leq \delta^{-1}$  
and $\| B^{-1}_{\delta }\| \leq  \| B^{-1} \|$,
 for the first part of the right hand side in \fr{hBdelta}, we obtain
 \begin{eqnarray}
 \| \hat{B}^{-1}_{\delta } -B^{-1}_{\delta } \| &\leq& 2 \| B^{-1}\|^2 \| \hat{B}-B\| \II (\bar{\Omega}_B) + 
\| \hat{B}^{-1}_{\delta }- B^{-1}_{\delta }\|  \II(\Omega_B) \nonumber \\
 & \leq & 2 \| B^{-1}\|^2  \| \hat{B}-B\| + 2\delta^{-1} \II(\Omega_B). \label{hBdB}
\end{eqnarray}
For the second part of \fr{hBdelta}  we derive $ \| B^{-1}_{\delta } - B^{-1} \| \leq \delta \| B^{-1} \|^2 $. 
Finally, combining  \fr{hadeltan},  \fr{hBdelta}  and \fr{hBdB}, we derive \fr{eq:hata_a}.
\\


 {\bf Proof of Lemma \ref{lem:randerror}.  } Recall that,  since the index set $K_{m,y}$ is finite,  by  \fr{variance}, one has
  \begin{eqnarray}  \label{R2}
R_{2n} &=& \EE\lkv   \sum_{k \in K_{m,y}} [(\h{a}_{\delta})_{m, k} - a_{m,k}] \varphi_{m,k} (y) \rkv^2  
=   O \lkr  2^m\ \EE \| \hat{a}_{\delta } - a \|^2  \rkr 
 \end{eqnarray}
In order to find an asymptotic upper bound for $\EE \| \hat{a}_{\delta } - a \|^2$, 
square the right-hand side in \fr{eq:hata_a} and find expectation in a view of 
Lemmas \ref{lem:large_dev}  and \ref{lem:hBhc}. Taking into account that it follows 
from \fr{Bexpan} and \fr{cexpan} that, as $m \rightarrow \infty$, one has  $\|c\| = O(2^{-m/2})$ and
$B^{-1}= (p(y))^{-1} I + O(2^{-m})$ where $p(y)>0$ is a non-asymptotic value, so that, 
$\|B^{-1} \| = O(1)$ and $\|B^{-1} \|^{-1} = O(1)$ as $m  \rightarrow \infty$, obtain
  \begin{eqnarray*} 
  \EE  \|  \hat{a}_{\delta } - a \|^2 & = & O \lkr \EE \| \hat{c}- c\|^2 + 2^{-m} [ 2^m n^{-1}  
+ \EE \| \hat{B}-B\|^2 + 2^{-m} n \PP(\Omega_B)] \right. \\
& + & 
\left. 2^m n^{-1}  \EE \| \hat{c}- c\|^2 +  \sqrt{\EE \| \hat{c}- c\|^4}  \sqrt{\EE \| \hat{B}-B\|^4} 
+ 2^{-m} n  \sqrt{\EE \| \hat{c}- c\|^4} \sqrt{\PP (\Omega_B)} \rkr
 \end{eqnarray*}
Applying \fr{large_devB} with $\tau = (2 \|B^{-1} \|^{-1} M \sqrt{\log n} )^{-1}  2^{-m/2} \sqrt{n} \rightarrow \infty$,
obtain that, as $n \rightarrow \infty$,
$$
 \PP(\Omega_B) = 0(n^{-h}) \ \ \ \mbox{for any} \ \ h>0, 
$$
i.e., $\PP(\Omega_B)$ tends to zero faster than any negative power of $n$. 
Then,  result of  Lemma \ref{lem:randerror} follows directly   from Lemma   \ref{lem:hBhc}  and 
the fact that  $2^m n^{-1} \rightarrow 0$ and $\|\gamma^{(2)} (m) \|^2\, 2^{2m} = o(n^3)$
as $n \rightarrow \infty$.\\


 {\bf Proof of Theorem \ref{th:lower_bounds}.  } 
Our goal is to construct a lower bound for the minimax risk when   $g \in  {\cal G}_r$.
In order to construct lower bound $\Delta_n (y)$ for \fr{rfun}, we use Theorem  2.7  of  
Tsybakov~(2008) which we reformulate here for the case of squared risk.

\begin{lemma} \label{Tsybakov} [Tsybakov~(2008), Theorem 2.7] Assume that $\Xi$ contains 
elements $\xi_0, \xi_1, \cdots, \xi_{n_0}$, $n_0 \geq 1$, such that \\
(i) $d(\xi_{i}, \xi_{i_0}) \geq 2 \chi >0$,    \  \  \  for $0 \leq i   < i_0 \leq n_0$;\\
(ii) $  P_{i} << P_ 0$, for $i = 1, \dots, n_0$, and
  \begin{eqnarray*}
K( P_{i}, P_0) \leq  C_{n_0},
  \end{eqnarray*} 
where $P_{i} = P_{\xi_{i}}$, $i = 0, 1, \dots, n_0,$ and $C_{n_0}$ is a positive constant. 
Then, for some absolute positive constant $C_K$, one has 
 \begin{eqnarray*}
 \inf_{\widehat{\xi}} \sup_{\xi \in \ \Xi} \EE_{\xi}  \left[d^2 (\widehat{\xi}, \xi)\right]  \geq C_K\, \chi^2.
 \end{eqnarray*}
\end{lemma}

Now,  consider  $\Xi= {\cal G}_r$, $d( f, g ) = \left| f(y) - g(y) \right|$, $n_0 =1$,
$p_1(x) = p_0(x) + \zeta k\left(h^{-1} (x-y) \right)$  and 
$\Psi_1(x)  = \Psi_0(x) + \zeta w_{h, y}(x)$.
%
Choose $\zeta= \zeta_0 \min\{h^r, \rho_r(h) \} = \zeta_0 h^{\max(r, r_1)}$
where $\zeta_0 \leq C_p$ is such that  $p_{1}(x) \geq 0$ and both $p_{1}(x)$ and $\Psi_{1}(x)$ are in  ${\cal G}_r$.
Let $t_i(y) = \Psi_i (y)/p_i(y)$, $i=0,1$.
Calculating the distance $d(t_{1}, t_0)$ at the fixed point $y$, due to $p_1(y) \geq p_0(y)/2$, we obtain
\begin{eqnarray*}
  d(t_{1}, t_0) & =&\left| \frac{ \Psi_0(y) + \zeta w_{h, y}(y)}{  p_0(y) + \zeta k(0)} - \frac{ \Psi_0(y)}{p_0(y)}\right|   
= \zeta \left|\frac{ w_{h, y}(y)p_0(y) - \Psi_0(y)k(0)}{p_0(y)[p_0(y) + \zeta k(0)]}\right|\\
  & \geq & \frac{\zeta}{2}\  \left|\frac{ w_{h, y}(y)p_0(y) - \Psi_0(y)k(0)}{p^2_0(y)}\right|.  
 \end{eqnarray*}
Hence,
$$
d(t_{1}, t_0) \geq  \lfi 
\begin{array}{ll}
     C \zeta, & \mbox{if}\ \ \ | w_{h, y}(y)| \leq C_0,\\
     C \zeta |w_{h, y}(y)|, & \mbox{if}\ \ \  \lim_{h \rightarrow 0}| w_{h, y}(y)| = \infty,
  \end{array} \right.
$$
so that  $\chi = C h^{\max(r, r_1) - r_2}$. In order to apply Lemma \ref {Tsybakov}, one needs to verify condition  (ii).  Observe that 
   \begin{eqnarray*}
   p_0(x_1, \cdots, x_n ) &=& \prod^{n}_{i=1}p_0(x_i), \quad
    p_{1}(x_1, \cdots, x_n ) = \prod^{n}_{i=1}\left[p_0(x_i) + \zeta k\left( \frac{x_i-y}{h}\right)\right].
     \end{eqnarray*}      
Then, due to  the fact that $\log ( 1 + x ) \leq x$, the Kullback divergence  between $p_{1}$ and $p_0$ is bounded  
\begin{eqnarray*}
K( p_{1}, p_0)   & =& \int \cdots \int \log \left\{\prod^{n}_{i=1} p_{1}(x_i)/p_0(x_i) \right\}\prod^{n}_{i=1}p_{1}(x_i)dx_i\\
& = & 
\sum^{n}_{i = 1} \int \log \left\{ \frac{ p_0(x) + \zeta k\left( \frac{x-y}{h}\right)}{p_0(x)} \right\} \left\{ p_0(x) + \zeta k\left( \frac{x-y}{h}\right) \right\}dx\\
& \leq & 
\sum^{n}_{i = 1} \int h\,  p_0(x)^{-1}\,   k\left( \frac{x-y}{h}\right) \lkv  p_0(x) + \zeta k\left( \frac{x-y}{h}\right) \rkv  dx\\
& = & 
n  \zeta^2 \int  p_0(x)^{-1}\,  k^2\left( \frac{x-y}{h}\right)  dx. 
\end{eqnarray*} 
 Now,  due to \fr{rhop}, one has  $n h  \zeta^2  h  \asymp  n h^{2\max(r, r_1) +1}$.
Therefore, $h \asymp n^{\frac{1}{2\max(r, r_1) +1}}$ and 
$$
\chi^2 = C n^{-\frac{2\max(r, r_1) -2r_2}{2 \max(r, r_1) +1}},
$$
which completes the proof.
\\


{\bf Proof of Theorem \ref{th:nonadaptive}.  } Validity of Theorem \ref{th:nonadaptive} follows 
directly from Lemmas~\ref{lem:syserror},~\ref{lem:randerror}.\\


\begin{lemma} \label{lem:error4} 
Let  $\delta^2 \sim n^{-1} 2^m$ and assumptions \fr{gamma_param_form} and \fr{uniform_cond} hold. 
Then,  under assumptions of Lemma \ref{lem:syserror}, as $n \rightarrow \infty$, 
\be \label{error4}
\EE\left|\widehat{t}_{{m}}(y)- t(y)  \right|^4 = O\left( n^{-2} 2^{2m} (  \gamma^2_m +1)^2+ 2^{-4mr} \right).
\ee
\end{lemma}

{\bf Proof of Lemma \ref{lem:error4}.  }
\begin{eqnarray}  \label{Et4}
\EE\left|\widehat{t}_{{m}}(y)- t(y)  \right|^4  
& \leq & 8\ \EE\left|\widehat{t}_{{m}}(y)- t_{{m}}(y)\right|^4 + 8 \left| t_ {{m}}(y)- t(y)  \right|^4. 
\end{eqnarray}
By Lemma \ref{lem:syserror},   as $m,n \rightarrow \infty$, one has  $\left| t_ {{m}}(y)- t(y)  \right|^4 =o(2^{-4{m}r})$.
For the first term in \fr{Et4}, note that 
$$
\EE\|\widehat{t}_{{m}}(y)- t_{{m}}(y) |^4 = O ( 2^{2m} \  \EE \|  \hat{a}_{\delta } - a \|^4).
$$ 
Upper bounds for $ \EE \|  \hat{a}_{\delta } - a \|^4$ can be derived using Lemma \ref{lem:adelta_dev}. 
Since, as $m,n \rightarrow \infty$,  $\|c\| = O(2^{-m/2})$ and $\|B^{-1} \| = O(1)$, 
it follows from  \fr{hBhc1}  and \fr{c_dev_moments}  that 
\begin{eqnarray*} 
 \EE \|  \hat{a}_{\delta } - a \|^4 & = & O \lkr \EE\| \hat{c}- c\|^4 \rkr + 
O \lkr 2^{-2m} \lkv \delta^4 + \EE\| \hat{B}-B\|^4 + \delta^{-4} P(\Omega_B) \rkv \rkr\\
& + & O \lkr \delta^4 \EE\| \hat{c}- c\|^4 + \sqrt{\EE\| \hat{c}- c\|^8} \sqrt{\EE\| \hat{B}-B\|^8} 
+ \delta^{-4}  \sqrt{\EE\| \hat{c}- c\|^8} \sqrt{P(\Omega_B)} \rkr\\
& = & O \lkr n^{-2}(\gamma_m^4 +1) \rkr =  O \lkr n^{-2}(\gamma_m^2 +1)^2 \rkr,
\end{eqnarray*} 
which completes the proof.
\\


{\bf Proof of Lemma \ref{lem:prob4}.  } 
 Denote  $R^2_{mn} =   \left[ \|\hat{B}_{\delta m}^{-1} \|^2 + \|\hat{B}_{\delta m_0}^{-1} \|^2 \right]^2  \rho^2_{mn}$ and 
 observe that
\begin{eqnarray*}
\PP( | \h{t}_m(y)-\h{t}_{m_0}(y) |\geq  \lambda R_{mn} )& \leq& 
\PP( | \h{t}_m(y)-{t}_{m}(y) | +  |{t}_m(y)-{t}(y) | \geq  0.5\, \lambda  R_{mn})\\
&+& \PP( | \h{t}_{m_0}(y)-{t}_{m_0}(y) | +  | {t}_{m_0}(y)-{t}(y) |\geq  0.5\, \lambda  R_{mn}).  
\end{eqnarray*} 
Since $m > m_0$ and $R_{mn}$ is an increasing function of $m$, one has $ |{t}_m(y)-{t}(y)| = o(2^{-mr})$ 
as $m\to\infty$  and $R_{mn} > R_{m_0 n}$. Therefore, it is sufficient to show that 
$$
\PP  \lkr |\h{t}_m(y)-{t}_{m}(y)|   \geq  0.5\, \lambda  R_{mn} - o(2^{-mr}) \rkr = O(n^{-2})
$$
 for any $m \geq m_0$. Taking into account that $|\h{t}_m(y)-{t}_{m}(y)|  \leq 2^{m/2} C_\varphi \|\hat{a}_m - a \|$ and     
$2^{-mr}/R_{mn} \rightarrow 0$ as $m,n \rightarrow \infty$, it is sufficient to show that 
\be \label{eq:to_show}
\PP  \lkr \|\hat{a}_m - a \|   \geq  2^{-m/2} (\lambda-1) R_{mn} /(2 C_\varphi) \rkr = O(n^{-2}),\ \ n \rightarrow \infty.
\ee

Recall that 
$\hat{B}^{-1}_{\delta}- B^{-1}= \hat{B}^{-1}_{\delta}(B-\hat{B}_{\delta})B^{-1}$, so that, for any $\delta>0$, one has 
$$
 \| \hat{B}^{-1}_{\delta}- B^{-1}\| \leq \|B^{-1}_{\delta}\|^2 \left( \| \hat{B}-B\| 
+ \delta_m \right) + 2 \| \hat{B}^{-1}_{\delta} \|^2 \| B^{-1}\| \left[ \| \hat{B} - B\|^2+ \delta^2_m \right].
$$
and, also,  
$$
 \| \hat{a}_{\delta}-a \| \leq \|\hat{B}^{-1}_{\delta} \|\| \hat{c}-c \| + \| \hat{B}^{-1}_{\delta} - B^{-1} \| \| c\|.
$$
Consequently, probability in \fr{eq:to_show} can be partition into three terms:
\be \label{eq:P123}
\PP  \lkr \|\hat{a}_m - a \|   \geq  2^{-m/2} (\lambda-1) R_{mn} /(2 C_\varphi) \rkr \leq \PP_1 + \PP_2 + \PP_3
\ee
where 
\begin{eqnarray*}
\PP_1 & = &   \PP \lkr \|\hat{B}^{-1}_{\delta} \|\| \hat{c}-c \| \geq \frac{\alpha_1 R_{mn} ( \lambda- 1)}{2^{m/2} 2C_{\varphi}} \rkr \\
\PP_2 & = &   \PP\lkr \|c\| [  \| \hat{B} - B\| + \delta_m] \geq \frac{ \alpha_2 
              \sqrt{1+ \gamma^2_m}\sqrt{\log n}( \lambda -1)}{2\sqrt{n} C_{\varphi}} \rkr\\
\PP_3 & = &   \PP\lkr \| B^{-1}\|  [  \| \hat{B} - B\|^2 + \delta^2_m] \geq \frac{ \alpha_3 
              \sqrt{1+ \gamma^2_m}\sqrt{\log n}( \lambda -1)}{4 \sqrt{n}C_{\varphi} \|c\|} \rkr\\
 \end{eqnarray*} 
and $\alpha_1, \alpha_2$ and $\alpha_3$ are positive constants such that $\alpha_1 + \alpha_2+ \alpha_3=1$. 

Applying \fr{large_devc} and taking into account that $\|\hat{B}_{\delta}\| \leq 2\|p\|_{\infty}$, obtain 
 \be \label{eq:P1}
 \PP_1   \leq     \PP  \lkr \| \h{c}-c\|^2 \geq \alpha_1 
\frac{ \sqrt{1+ \gamma^2_m}\sqrt{\log n}( \lambda -1)}{4  \|p\|_{\infty}\sqrt{n} C_{\varphi}} \rkr  \leq  2M n^{-\tau_1} 
\ee
where $\tau_1= (128 M  C^2_{\varphi}\|p\|^3_{\infty} M)^{-1}\, \alpha^2_1(\lambda-1)^2$.
 Recalling that $\|c\| \leq 2M \| \Psi\|_{\infty} 2^{-m/2}$, using formula  \fr{large_devB} and 
taking into account that $1 - 4M  \|\Psi\|_{\infty} C_{\varphi}/(\alpha_2 (\lambda -1) > 1 - \nu_1$
for any small positive constant $\nu_1$ as $n \rightarrow \infty$,    we derive
 \begin{eqnarray}
 \PP_2 & \leq & \PP \lkr  \| \hat{B} - B\| \geq   \frac{ ( \alpha_2 \lambda -1) 2^{m/2}\, \sqrt{1+ \gamma^2_m}\sqrt{\log n}}
{4M  \|\Psi\|_{\infty} C_{\varphi}\, \sqrt{n} }- \frac{2^{m/2}}{\sqrt{n}} \rkr \nonumber\\
 & \leq&  \PP \lkr  \| \hat{B} - B\| \geq \frac{ M 2^{m/2}  \sqrt{\log n}}{\sqrt{n}}\ \frac{ \alpha _2 (\lambda -1) \sqrt{ 1-\nu_1}}{4M^2} \rkr
\leq 2M^2 n^{-\tau_2} \label{eq:P2}
 \end{eqnarray} 
where 
$\tau_2= (128 M^4 C_{\varphi}^2 \|\Psi\|_{\infty}^2  \|\varphi\|_\infty^2 \|p\|_{\infty})^{-1}\, \alpha^2_2(1-\delta_1)(1-\lambda)^2$.
 
In order to find an upper bound for $\PP_3$, recall that  $\|B^{-1}\| \leq   2M/p(y)$ and $\|c\| \leq 2M \| \Psi\|_{\infty} 2^{-m/2}$.
Also, note that $p(y) \geq (\log n)^{-1/2}$, for any fixed $y$,  as $n \rightarrow \infty$.
Therefore, applying \fr{large_devB} and taking into account that, due to \fr{eq:mn} and $m \leq m_n$, one has 
$ (2^{-m/2} \sqrt{n} -1)/\log n > 1 - \nu_2$ for any small positive constant $\nu_2$ as $n \rightarrow \infty$,  derive
 \begin{eqnarray}
\PP_3 & \leq & \PP \lkr  \| \hat{B} - B\|^2   \geq \frac{ \alpha_3 ( \lambda -1) \sqrt{1+ \gamma^2_m}\sqrt{\log n}}
{2 \sqrt{n}C_{\varphi}\| B^{-1}\|  \|c\|} - 2^m n^{-1} \rkr \nonumber\\
& \leq &  
\PP \lkr \| \hat{B} - B\|^2\leq  \frac{2^m M^2 \log n}{n} \frac{ \alpha_3( \lambda -1)}{16 C_{\varphi} M^3  \|\Psi\|_{\infty}} \rkr
\leq 2M^2 n^{-\tau_3} \label{eq:P3}
 \end{eqnarray}   
where $\tau_3=  (128 M^3 C_{\varphi} \|\Psi\|_{\infty}   \|\varphi\|_\infty^2 \|p\|_{\infty}^2)^{-1}\, \alpha_3(\lambda-1)(1-\nu_2) $.

Now, in order to complete the proof, combine \fr{eq:P123} -- \fr{eq:P3} and choose $\alpha_i$, $i=1,2,3,$ 
such that $\tau_i \geq 2$ for  $i=1,2,3,$ and $\PP_1+ \PP_2 + \PP_3$ takes minimal value. 
\\


{\bf Proof of Theorem  \ref{th:adaptive}.} 
First, let us show that $\EE \lkv \|\hat{B}_{\delta m}^{-1} \|^2 + \|\hat{B}_{\delta m_0}^{-1} \|^2 \rkv^2 = O(1)$ as $m,n \rightarrow \infty$,
so that asymptotic relation \fr{eq:Delta1_portion} holds. Indeed, for $m_1 \leq m \leq m_0$ and any fixed $y$, one has
\begin{eqnarray*}
\|\hat{B}_{\delta m}^{-1} \| & \leq & \|\hat{B}_{\delta m}^{-1} - B_{\delta m}^{-1} \| + \| B_{\delta m}^{-1} \|  \\
& \leq &  2  \| B_{\delta m}^{-1} \|^2  \|\hat{B}_{\delta m}  - B_{\delta m}  \| + 2 \delta_m^{-1} \II(\Omega_m) + \| B_m^{-1} \| 
\end{eqnarray*}
where $\Omega_m$ is defined in Lemma \ref{lem:adelta_dev}. Then, 
$$
\EE \|\hat{B}_{\delta m}^{-1} \|^4 = O \lkr \EE  \|\hat{B}_{\delta m}  - B_{\delta m}  \|^4 + \delta_m^{-4} \PP(\Omega_m) 
+ \| B_m^{-1} \|^4 \rkr = O(1),
$$
so that both \fr{eq:Delta1_portion} and \fr{Delta1} are valid. 
 
In order to find an upper bound for $\Delta_2$, note that by Lemmas \ref{lem:error4} and \ref{lem:prob4}, one has 
\begin{eqnarray*} 
\Delta_2 & = & \EE[\left|\widehat{t}_{\widehat{m}}(y) - t(y)  \right|^2 \II( \widehat{m} > m_{0} )]   \\
&\leq & 
\sum^{m_n}_{l= m_0+1} \sqrt{  \EE\left|\widehat{t}_{\widehat{m}}(y)
- t(y)  \right|^4}\, \sqrt{ \PP( | \h{t}_j(y)-\h{t}_{m_0}(y) |^2\geq \lambda^2\rho^2_{jn}\left[ \|\hat{B}_{\delta j}^{-1} \|^2 + \|\hat{B}_{\delta m_0}^{-1} \|^2 \right]^2)}  
 \end{eqnarray*}
so that $\Delta_2 = O(n^{-1}) = o\lkr \rho^2_{m_0n}  \rkr$ as $n \rightarrow \infty$.



\begin{thebibliography}{99}

\bibitem{brown1}
  Brown, L. D., Greenshtein, E. (2009). 
Nonparametric empirical Bayes and compound decision approaches to estimation 
of a high-dimensional vector of normal means. {\it Ann. Statist.}, {\bf 37}, no. 4, 1685–1704.


\bibitem{brown2}
 Brown, L., Gans, N., Mandelbaum, A., Sakov, A., Shen, H., Zeltyn, S., Zhao, L. (2005). 
Statistical analysis of a telephone call center: a queueing-science perspective. 
{\it J. Amer. Statist. Assoc.}, {\bf100} , no. 469, 36–50, 


\bibitem{carlin}
Carlin, B.P., Louis, T.A.(2000)
{\it Bayes and empirical Bayes methods for data
   analysis}, Second ed., Chapman\& Hall/CRC, Boca Raton.
   
   \bibitem{Sasel1}
Casella, G. (1985)
An introduction to empirical Bayes data analysis. {\it The American Statistician}, {\bf 39},   83--87. 


\bibitem{datta1}
Datta, S. (1991)
Nonparametric empirical Bayes estimation with $O(n\sp {-1/2})$ rate of
   a truncation parameter. {\it Statist. Decisions}, {\bf 9},   45--61. 


\bibitem{datta}
Datta, S.  (2000)
Empirical Bayes estimation with non-identical components. 
{\it  J. Nonparametr. Statist.}, {\bf 12},  709--725. 

\bibitem{EfMor}
Efron, B., Morris, C. N. (1977)
Stein's paradox in statistics. {\it Scientific American}, {\bf 236},   119--127. 


\bibitem{ghosh1}
Ghosh, M., Lahiri, P. (1987) 
Robust empirical Bayes estimation of means from 
stratified samples. {\it J. Amer. Statist. Assoc.}, 
{\bf 82}, 1153--1162.


\bibitem{ghosh2}
Ghosh, M., Meeden, G. (1986)
Empirical Bayes estimation in finite population sampling.
{\it J. Amer. Statist. Assoc.}, {\bf 81}, 1058--1062.


\bibitem{huang}
Huang, S.Y. (1997) 
Wavelet based empirical Bayes estimation for the uniform distribution.
{\it   Statist. Probab. Lett.}, {\bf 32},   141--146.
 

\bibitem{jiang}
Jiang, W.,   Zhang, C.-H. (2009)
General maximum likelihood empirical Bayes estimation of normal means.
{\it Ann. Statist.},  {\bf 37},   1647-1684. 


\bibitem{karunamuni}
Karunamuni, R.J. (1996)
Optimal rates of convergence of empirical Bayes tests for the continuous one-parameter exponential family.
{\it Ann. Statist.},  {\bf  24},   212-231.


\bibitem{karusingh}
Karunamuni, R. J., Singh, R. S., Zhang, S. (2002) 
On empirical Bayes estimation in the location family.
{\it J. Nonparam. Stat. }, {\bf 14},  435--448. 


\bibitem{karuzhang}
Karunamuni, R.J.; Zhang, S. (2003)
Optimal linear Bayes and empirical Bayes estimation and
prediction of the finite population mean. 
{\it J. Statist. Plan. Inference}, {\bf 113},   505--525.


\bibitem{lep1} 
Lepski, O.V. (1991). 
Asymptotic mimimax adaptive estimation. I: Upper bounds. Optimally adaptive estimates.
{\it Theory Probab. Appl.}, {\bf 36}, 654-659.

\bibitem{lep2}  
Lepski, O. V., Mammen, E., Spokoiny, V. G. (1997). 
Optimal spatial adaptation to inhomogeneous smoothness: An approach based on 
kernel estimators with variable bandwidth selectors.  
{\it Ann.\ Statist}, {\bf 25}, 929-947.


\bibitem{li}
Li, J., Gupta, S.S., Liese, F.(2005).
Convergence rates of empirical Bayes estimation in exponential family.
{\it J. Statist. Plan. Inference}, {\bf 131}, 101-115.
 

\bibitem{liang}
Liang, T. (2000)
On an empirical Bayes test for a normal mean.
{\it Ann.\ Statist}, {\bf 28},    648-655.

 

\bibitem{ Louis}
Louis, T. A. (1984)
Estimating a population of parameter values using Bayes and empirical Bayes methods.
{\it J. Amer. Statist. Assoc},  {\bf 79},   393-398. 




\bibitem{ma}
Ma, Y., Balakrishnan, N. (2000)
Empirical Bayes estimation for truncation parameters. 
{\it J. Statist. Plan. Inference}, {\bf 84}, 111--120.
 


\bibitem{mallat1999}
Mallat, S. (1999). {\em A Wavelet Tour of Signal Processing}. 
2nd Edition, Academic Press, San Diego.

\bibitem{maritz}
Maritz, J.S., Lwin, T. (1989). 
{\it Empirical Bayes Methods} (2nd ed.), 
Chapman \& Hall, London.

\bibitem{Morris}
Morris, C. N. (1983)
Parametric empirical Bayes inference: Theory and applications.
{\it J. Amer. Statist. Assoc.},  {\bf 78},   47-65. 



\bibitem{nog} 
Nogami, Y. (1988)  
Convergence rates for empirical
Bayes estimation in the uniform $U(0,\theta)$ distribution.
{\it Ann. Stat.}, {\bf 16}, 1335-1341.


\bibitem{penskaya}
Penskaya, M. (1995) On  the lower bounds for mean square error of
empirical Bayes estimators.  {\it Journ. of  Math. Sciences},
{\bf 75},  1524--1535.

\bibitem{pen1}
Pensky, M. (1997).  
A general approach to nonparametric empirical Bayes estimation.  
{\it Statistics}, {\bf 29}, 61--80.


\bibitem{pen2}
Pensky, M. (1997). Empirical Bayes estimation of a location parameter.
{\it Statist. Decisions}, {\bf 15}, 1--16.

\bibitem{pen3}
Pensky, M. (1998) Empirical Bayes estimation based on wavelets.
{\it Sankhy\={a}}, {\bf  A60}, 214--231.

 

\bibitem{pen6}
Pensky, M. (2000)  
Adaptive wavelet empirical Bayes estimation of a location or a scale parameter. 
 {\it J. Statist. Plan. Inference}, {\bf 90}, 275 --292.

 
 
\bibitem{pen7}  Pensky, M. (2002). 
Locally adaptive wavelet Bayes estimation of location parameter. 
{\it Ann. Inst. Statist.  Math}, {\bf 54}, 83-99. 



\bibitem{pen8}
Pensky, M. (2003)
Rates of convergence of   empirical Bayes tests for a normal mean.
{\it Journal of Statistical Planning and Inference}, {\bf 11}, 181--196.


\bibitem{penal} 
Pensky, M., Alotaibi, M.  (2005)
Generalization of linear empirical Bayes estimation via wavelet series.
{\it Statistics and Decisions},  {\bf 23}, 181--198.


\bibitem{pen5}
Pensky,M., Ni, P. (2000) 
Extended linear empirical Bayes estimation.
{\it Communications in Statistics - Theory and Methods}, {\bf 29}, 579 -- 592.


 \bibitem{raykar} Raykar, V.,  Zhao,  L. (2011).
Empirical Bayesian thresholding for sparse signals using mixture loss functions
{\it Statistica Sinica}, {\bf 21},  449-474.


\bibitem{rob1}
Robbins, H. (1955)  
An empirical Bayes approach to statistics. 
In {\it Proceedings of the Third Berkeley Symposium on Mathematical 
Statistics and Probability}, (Vol. 1), University of California Press, Berkeley, 157--163.


\bibitem{rob2}
Robbins, H. (1964)  
An empirical Bayes approach to statistical decision problems. 
{\it Ann. Math. Statist.}, {\bf 35}, 1--19.



\bibitem{rob3}
Robbins, H. (1983)   
Some Thoughts on Empirical Bayes Estimation.  
{\it Ann. Statistics}, {\bf 11}, 713--723.


\bibitem{singh}
 Singh, R. S. (1976)
Empirical Bayes estimation with convergence rates in noncontinuous
   Lebesgue exponential families. {\it Ann. Statist.}, {\bf 4},   431--439.


\bibitem{snh}  
Singh,~R.S. (1979)  
Empirical Bayes estimation in  Lebesgue
-exponential families with rates near the best possible rate.
{\it Ann. Statist.}, {\bf 7}, 890--902. 


  \bibitem{tsybakov}
Tsybakov, A.B. (2008)
{\it Introduction to Nonparametric Estimation}, Springer, New York.


\bibitem{walter}
Walter, G.G., Hamedani, G.G. (1991)  
Bayes empirical Bayes estimation for natural exponential families 
with quadratic variance function.  
{\it Ann. Statist.}, {\bf 19}, 1191--1224.


\bibitem{walshen} 
Walter, G.G., Shen, X. (2001) 
{\it Wavelets and Other Orthogonal Systems}. 
Chapman and Hall/CRC, Boca Raton.


\end{thebibliography}
  \end{document}